\documentclass[11pt]{article}    
\usepackage{amsmath,amssymb,latexsym,epsfig,subfigure,multirow,longtable}  
\begin{document}

\title{On the validity of ``A proof that the discrete singular convolution
(DSC)/Lagrange-distributed approximation function (LDAF) method is inferior
to high order finite differences''}
 
\author{ G.~W.~Wei$^{1,2}$\footnote{ 
Corresponding author. Tel: (517)3534689, Fax: (517)4321562, 
Email: wei@math.msu.edu} ~
and Shan~Zhao$^1$ 
\\
%\address{
\small \tt
$^1$Department of Mathematics,\\ \small \tt
   Michigan State University, 
   East Lansing, MI 48824, 
      USA\\ \small \tt
$^2$Department of Electrical and Computer Engineering,\\ \small \tt
    Michigan State University,  East Lansing, MI 48824, USA \\ \small \tt
} 
 %\date{\empty} 
 
\maketitle 
\begin{abstract} 
 A few families of counterexamples are provided
to ``A proof that the discrete singular convolution
(DSC)/Lagrange-distributed approximation function (LDAF) method is inferior
to high order finite differences'', Journal of Computational Physics, {\bf214}, 538-549 (2006).
\end{abstract}

%\vspace*{1cm} 

%%%%%%%%%%%%%%%%%%%%%%%%%%%%% 
%%% Start of Introduction %%% 
%%%%%%%%%%%%%%%%%%%%%%%%%%%%% 

\section{Introduction}
Recently, Boyd  published a  paper entitled ``A proof that the discrete singular convolution 
(DSC)/Lagrange-distributed approximating function (LDAF) method is inferior to high order finite 
differences'' \cite{boyd06}, which will be referred to as ``Proof'' throughout the present work. 
The sole purpose of the present paper is to analyze the validity of many  of the statement and claims in
``Proof''. To this end, we provide a wide variety of counterexamples. It is pointed out that we have no 
intention to claim the superiority of the DSC algorithm. In the rest of the Introduction, 
we shall  set the scope of our discussions.

For the DSC algorithm, only a special DSC kernel will be discussed. ``Proof'' 
makes very general claims about the DSC algorithm. Unfortunately, it ignores the fact that  
the DSC algorithm can be realized by many different kernels that might behave very differently 
from each other \cite{weijcp99}. What is discussed in ``Proof'' is a special DSC kernel,
the regularized Shannon kernel (RSK). We therefore limit our discussion to the DSC-RSK method  
exclusively in the present work.

One of ``Proof'''s central claims is about spectrally-weighted differences. ``Proof'' states 
that ``Although we shall not perform detailed comparisons between DSC and spectrally-weighted 
differences, the good performance of DSC for high $k$ and small $a$ is an accident. It seems 
likely that for $f(x)$ which are known to have spectra concentrated between $K = \pi/2$ and 
$K = \pi$, one could obtain higher accuracy from spectrally-weighted differences than from DSC''.
It is seen that ``Proof''  does not provide any detailed comparisons between the DSC and 
spectrally-weighted differences. Therefore,  ``Proof''  has to present its claim about the
superiority of the spectrally-weighted differences over the DSC  as a shaky
speculation: ``It seems likely that ...''.  Unfortunately, this unfounded superiority was 
a part of ``Proof'''s main claims. Although there was no proof to begin with, ``Proof''  still 
concludes that ``DSC/LDAF methods are {\it never} the best way to approximate 
derivatives on a stencil of a given width'' in its abstract and  ``the DSC/LDAF method is 
{\it never} the method of choice for approximating derivatives'' in its summary. 
In this work, we will analyze the validity of this claim in detail. Since ``Proof'''s   
claim about spectrally-weighted differences is {\it very general}, logically, it is 
{\it sufficient} for us to discredit this general claim by considering a special case. 
To this end, we choose Boyd's sech function based spectrally-weighted difference (Sech) 
\cite{Boyd94}, for which a construction procedure  was emphasized and outlined in ``Proof''. 
It is pointed out that our results in this work carry no implication on the performance of any 
other spectrally-weighted differences.

The other part of ``Proof'''s central claims is about the standard high order 
finite difference (FD) scheme. ``Proof'' classifies functions into two classes, namely, those 
with small (Fourier) wavenumbers and those with large (Fourier) wavenumbers, and then 
argues that the standard high order finite difference (FD) scheme is more accurate than the DSC 
in differentiating functions with considerable amplitude in small wavenumbers, while spectrally-weighted 
differences, e.g., the Sech method,  are superior for functions with (solely) large wavenumbers. 
``Proof'' claims that ``Consequently, DSC/LDAF methods are {\it never} the best way to approximate 
derivatives on a stencil of given width''. 

It should be  pointed out that no part of ``Proof'''s claims is based on rigorous error analysis.
Error analysis for the DSC-RSK method was given earlier by Qian \cite{Qian}, and was not mentioned 
in ``Proof''. Although error expressions for both the FD and DSC-RSK methods are presented, 
these expressions are not compared and do not directly support ``Proof'''s  claims. Instead, 
``Proof'''s claim about the superiority of the FD method over the DSC-RSK method is based on numerical experiments 
with a limited set of parameters, and on some informal arguments. In this work, we present a counterexample in 
differentiating $e^{ikx}$ with a small wavenmuber $k$, see Fig. \ref{fig.dtest}. In differentiating 
a function with exponentially decaying amplitudes in wavenumbers $k$, the DSC-RSK method is more accurate than 
the FD method over all the stencils examined, and is up to six orders more accurate than the FD method and the 
Sech method at some large stencils,  see Fig. \ref{fig.ftest2}. For problems with 
considerable amplitudes in a wide range of wavenumbers as shown in Fig. \ref{fig.Helmsolu} (b), the 
DSC-RSK method outperforms the FD method up to a factor of $10^{10}$, see Table \ref{table.Helm2}. In differentiating 
$e^{ikx}$ with a relatively large wavenmuber, errors in Boyd's spectrally weighted difference, the Sech method, 
are up to $10^8$ times larger than those of the DSC-RSK method, see Fig. \ref{fig.dtest} (b). 
In a variety of other counterexamples, the DSC-RSK method outperforms the FD and Sech methods up to factors of multiple 
orders of magnitude for these problems.

Apart from the Sech and  FD methods, ``Proof'' places great emphasis about a few other generalized finite 
differences, namely, 
Boyd's finite difference (Boyd's FD) \cite{Boyd94},
Boyd's Euler-accelerated sinc algorithm (Euler) \cite{Boyd91,Boyd94}, and
Boyd's modified Euler-accelerated sinc algorithm (MEuler) \cite{Boyd91,Boyd94}. 
Some detailed expressions and/or numerical procedures have been given in ``Proof'' for
these methods. Moreover, some of ``Proof'''s claims involve these methods.
Therefore, we will perform detailed comparisons of these methods with the 
DSC-RSK method as well. These comparisons will enhance our understanding of high order methods, 
as well as the validity of ``Proof'''s claims.

This paper will focus on the DSC algorithm and directly dismiss ``Proof'''s claim about LDAF.
In its title and abstract, ``Proof'' refers o LDAF as ``Lagrange-distributed approximating 
function'', while in its introduction,  ``Proof'' refers to LDAF as  ``linear distributed approximating 
functional''. The term ``linear distributed approximating functional'' was attributed to 
Hoffman et al. \cite{Hoffman}. Unfortunately, ``linear distributed approximating functional (LDAF)''
simply does not exist \cite{Hoffman}.  The paper by Hoffman et al. \cite{Hoffman} concerned
the so called  ``distributed 
approximating functional'' (DAF). Apart from this confusion, no  single expression, 
no detailed analysis, nor even any correct literature reference was given to LDAF throughout ``Proof''. 
Thus, we have no clear idea what ``Proof''  has really proved regarding LDAF. For these reasons, 
we have to ignore ``Proof'''s  claim about LDAF in the rest of this paper.

The organization of this paper is the follows. In Section \ref{analysis},
we analyze  ``Proof'''s major claims in detail, while the  
analysis of ``Proof'''s other  claims will be given in Section \ref{further}.
Relevant methods and algorithms are briefly described in Section \ref{methods}. 
Section \ref{results} is devoted to numerical results and 
discussions. Finally some concluding remarks are given.

\section{Analysis of ``Proof'''s central claims}\label{analysis}

For convenience, we quote the abstract of ``Proof'', which contains all major claims.

\begin{itemize}

\item {\sc [``Proof''- Abstract]}
``Finite differences approximate the $m$th derivative of a function u(x) by a series 
$\sum_{j=-M}^M d_j^{(m)}u(x_j)$,  where the $x_j$ are the grid points. The closely-related 
discrete singular convolution (DSC) and Lagrange-distributed approximating function (LDAF) 
methods, treated here as a single algorithm, approximate derivatives in the same way as 
finite differences but with different numerical weights that depend upon a free parameter $a$.
By means of Fourier analysis and error theorems, we show that the DSC is worse than the standard 
finite differences in differentiating $\exp^{(ikx)}$ for all $k$ when $a\geq a_{\rm FD}$ where  
$a_{\rm FD}=1/\sqrt{M+1}$ with $M$  as the stencil width is the value of the DSC parameter that makes 
its weights most closely resemble those of finite differences.  For $a < a_{\rm  FD}$, the DSC errors 
are less than finite differences for $k$ near the aliasing limit, but much, much worse for smaller $k$. 
Except for the very unusual case of low-pass filtered functions, that is, functions with negligible 
amplitude in small wavenumbers $k$, the DSC/LDAF is less accurate than finite differences for all 
stencil widths $M$. So-called ``spectrally-weighted'' or ``frequency-optimized'' differences are superior 
for this special case. Consequently, DSC/LDAF methods are {\it never} the best way to approximate 
derivatives on a stencil of a given width.'' 

{\sc [Analysis]}
First, we note that Fourier analysis is the sole technique used in ``Proof''. 
It is  improper to draw conclusions regarding the
performance  of a numerical method for solving partial differential equations solely from 
its performance in differentiating $e^{ikx}$. We will elaborate on these points in Section \ref{further}.

Second, ``Proof'' claims that ``we show that the DSC is worse than the standard 
finite differences in differentiating $\exp^{(ikx)}$ for all $k$ when $a\geq a_{\rm FD}$ where  
$a_{\rm FD}=1/\sqrt{M+1}$ with $M$  as the stencil width is the value of the DSC parameter that makes 
its weights most closely resemble those of finite differences''.
There is clearly no need to refute such a claim.
Because $a$ is a free parameter, no practitioner would impose the 
condition $a\geq a_{\rm FD}$, but instead, an optimal choice for this parameter is sought and typically 
this is a value $a <  a_{\rm FD}$

Third,  although here ``Proof''  avoids specifying what is meant by ``small $k$'' and what is the remainder,
it does give two intervals $|K|<\frac{\pi}{2}$ and $|K|>\frac{\pi}{2}$ in Section 6, where
$K=kh$, and $h$ is the grid spacing. In Fig. \ref{fig.dtest} (a), counterexamples are given to show that the 
DSC method outperforms the FD method at a small wavenumber, $K=\frac{9}{20}\pi$, over a wide range of stencils.

Fourth, ``Proof''  claims that  ``For $a < a_{\rm  FD}$, the DSC errors 
are less than finite differences for $k$ near the aliasing limit, but much, much worse for smaller $k$. 
Except for the very unusual case of low-pass filtered functions, that is, functions with negligible 
amplitude in small wavenumbers $k$, the DSC/LDAF is less accurate than finite differences for all 
stencil widths $M$''. This claim offers great details about the use of parameter $a$, 
the behavior with regarding to wavenumbers $k$, the prescription of amplitude, and the behavior with respect to stencil 
widths. One might expect to find some verification of these statements in ``Proof''. 
Unfortunately, no verification is given.  Indeed none is possible for for such a general statement.
We have designed problems involving a wide 
range of wavenumbers, in particular, including cases whose amplitudes in small wavenumbers $k$ are either
exponentially larger, or at least as large as those in large wavenumbers $k$, see Figs. \ref{fig.ftest1}(b), 
\ref{fig.BVPsolu3} (b), \ref{fig.Helmsolu} (b) and \ref{fig.hyperb2} (a). Our results for these four problems 
are given in Fig. \ref{fig.ftest2}, Fig. \ref{fig.bvp3}, Table \ref{table.Helm2}, and Fig. \ref{fig.hyperb2} (b), 
respectively. The DSC-RSK method outperforms the FD method up to $10^6$, $10^{5}$, $10^{10}$ and $10^{6}$ times, respectively
for these four problems. In  all four examples, the DSC outperforms the FD method over all stencil 
widths examined. We note that ``Proof''  restricted its analysis to only differentiating functions, 
while our last three counterexamples concern the solution of  boundary value, eigenvalue and 
unsteady hyperbolic problems. 
   
Fifth, for large wavenumbers, while admitting the superiority of the DSC-RSK method over the standard FD scheme, 
``Proof''  insists that spectrally weighted differences are superior to the DSC-RSK scheme. 
To counter this, we just need to analyze a special case --- the Sech scheme.
In Fig. \ref{fig.dtest} (b), we show that for differentiating $e^{ikx}$ at $K=\frac{3}{5}\pi$, 
which is a function with no amplitude in small wavenumbers $k$, the DSC-RSK method is up to $10^8$ times 
more accurate than Boyd's Sech method. Another example that involves solely large wavenumbers is given in
Fig. \ref{fig.BVPsolu2}. Errors in the Sech scheme are up to a factor of $10^{5}$ times larger than 
those of the DSC-RSK method, see Fig. \ref{fig.bvp2}. In fact, in this case, the DSC-RSK method outperforms the 
Sech scheme over all stencil widths examined.

In fact, 
a large number of counterexamples exist. In Section \ref{results}, we show that up to a factor of $10^{6}$
(see Tables \ref{table.Helm1} and \ref{table.Helm2}, Figs. \ref{fig.hyperb1} (b) and 
\ref{fig.hyperb2} (b)), the DSC-RSK method outperforms the FD method and Boyd's four methods in various 
problems examined.

\end{itemize}

\section{Methods and algorithms}\label{methods}

To avoid any confusion, this section presents a brief description of the numerical
methods that are studied in this work. These methods were chosen as the subject
of detailed discussions in ``Proof''  \cite{boyd06}, and can be expressed
in the following form
\begin{equation}\label{center}
\frac{d^n}{d x^n} u(x) \approx 
\sum_{j=-M}^{M} \delta_{j}^{(n)} u(x+x_j), 
\end{equation}
where $(2M+1)$ is the stencil width, $x_j=jh$ with $h$ being the spacing, and 
$\delta_{j}^{(n)}$ is the differencing kernel for $n$th order derivative
approximation. Note that Eq. (\ref{center}) gives an interpolation approximation
when $n=0$, while in general, the differentiation point $x$ may not be on-grid,
i.e., $x \ne x_j$ for $j=-M,\ldots,M$.

Sum-accelerated pseudospectral methods are based on the Sinc kernel
\begin{equation}\label{sinc}
d^{\rm (1),~sinc}_{j} = \left\{  
\begin{array}{ll}
(-1)^{j+1}/hj, & j=1,2,\ldots \\
0,            & j=0,
\end{array} \right. \quad \quad
d^{\rm (2),~sinc}_{j} = \left\{  
\begin{array}{ll}
2(-1)^j/h^2j^2, & j=1,2,\ldots \\
-\pi^2/3h^2,            & j=0.
\end{array} \right.
\end{equation}
Various accelerations have been proposed by Boyd for on-grid (i.e., differentiation
only at $x=0$) derivative approximation \cite{Boyd91,Boyd94}
\begin{equation}
\frac{d^n}{d x^n} u(0) \approx 
\sum_{j=-M}^{M} w_{Mj} d^{\rm (n),~sinc}_{j} u(x_j), \quad n=1,2,
\end{equation}
where $w_{Mj}$, $j=-M,\ldots,M$, are the acceleration weights.

Notation of various higher order approaches is given as the follows.
\begin{itemize}
\item FD: the standard finite difference method. 
The FD weights are calculated via Fornberg's code \cite{Fornberg}.

\item Boyd's FD: the finite difference method
generated by applying sum acceleration, with a certain set of
acceleration weights, to the standard sinc pseudospectral method
\cite{Boyd94}. The acceleration weights are given as
\begin{equation}
w_{Mj}=\left\{ 
\begin{array}{ll}
\frac{(M!)^2}{(M-j)! (M+j)!}, & j=1,2,\ldots,M, \\
\frac{6}{\pi} \sum_{k=1}^{M} \frac{1}{k^2},        & j=0.
\end{array} \right.
\end{equation}
The factorial involved in sum-accelerated differences is calculated
by calling the ``dgamma'' function of the SPECFUN package
\cite{Cody}.

\item Euler: the Euler-accelerated sinc algorithm
\cite{Boyd91,Boyd94}.
The acceleration weights are given as
\begin{eqnarray}
w_{Mj} & = & \sum_{k=j}^M \mu_{Mk} \label{Euler} \nonumber \\
\mu_{Mk} & = & \frac{M!}{2^n k! (M-k)!}, \quad k=0,1,\ldots, M. 
\end{eqnarray}

\item MEuler: the modified Euler-accelerated sinc algorithm
\cite{Boyd91,Boyd94}.
The acceleration weights are the same as those of Euler scheme,
except that $w_{M0}$ is modified so as to balance the contribution of all the other
weights in the second order derivative approximation.
\begin{equation}
w_{M0}^{\rm (mod)} = \frac{12}{\pi^2} \sum_{j=1}^{M} 
\frac{(-1)^{j+1} w_{Mj}}{j^2}.
\end{equation}

\item Sech: the spectrally weighted least squares differences with
the sech weights \cite{Boyd94}.
A difference scheme is said to be ``spectrally-weighted'' if 
the kernel $\delta_{j}^{(n)}$ of Eq. (\ref{center}) is chosen to minimize
\begin{equation}
\int_{-\pi}^{\pi} \omega(K) \left| K^n - 
\sum_{j=-M}^M \delta_{j}^{(n)} \exp( i j K) \right|^2 dK,
\end{equation}
where $K$ is the scaled wavenumber, $K=kh$. 
To calculate the weights, one first lets
\begin{equation}
\phi_j(K)=\left\{ \begin{array}{ll}
\sin (jK), & n~\textrm{odd}, \\
\cos (jK), & n~\textrm{even}.
\end{array} \right.
\end{equation}
Define the matrix and vector elements
\begin{equation}
G_{ij}= \int_{-\pi}^{\pi} \phi_i(K) \phi_j(K) \omega(K) dK, \quad
\chi_i = \int_{-\pi}^{\pi} \phi_i(K) K^n  \omega(K) dK.
\end{equation}
In the present study, these integrals are calculated up to double precision
limit by using Gaussian quadrature with 512 nodes. The approximation
accuracy is confirmed by benchmarking with the extended Simpson's rule
with a very large number of nodes.
The differentiation kernel is finally obtained by solving the linear system
\begin{equation}
\vec{\vec{G}} \vec{\delta} = \vec{\chi}.
\end{equation}
Following Boyd's suggestion \cite{Boyd94}, 
this least squares algebraic system is solved by using
the standard singular value decomposition (SVD) solver \cite{Press}.

If the Fourier spectrum of the function under study is known, one can
simply choose the least squares weight function $\omega(K)$ as the Fourier spectrum,
to attain an ideal approximation. However, such an information is usually
unavailable in advance. Thus, the best one can do is making use of some empirical
spectrum information to optimize $\omega(K)$.
For example, for a smooth function whose spectrum decays exponentially with $K$ for
sufficiently large $|K|$, Boyd introduced the sech weights \cite{Boyd94} 
\begin{equation}
\omega (K) = \textrm{sech} (K \pi / 2 D),
\end{equation}
where $D$ is an adjustable constant.

\item Sinc: the truncated Sinc pseudospectral method, i.e., 
the difference kernel $\delta_{j}^{(n)}$ of Eq. (\ref{center}) is chosen as
$d^{\rm (n),~sinc}_{j}$ given in Eq. (\ref{sinc}).

\item DSC-RSK: the regularized Shannon kernel (RSK) of the discrete singular
convolution (DSC) algorithm \cite{weijcp99}.
Although many other DSC kernels can be similarly employed,
regularized Shannon kernel (RSK) \cite{weijcp99} is employed in  
the present study, 
\begin{equation} 
\delta(x-x_k) = \frac{\sin \frac{\pi}{h} 
(x-x_k)}{\frac{\pi}{h} (x-x_k)} e^{-\frac{(x-x_k)^2} 
{2 \sigma^2}}, \quad k=-M,\ldots,M,
\end{equation} 
where the parameter $\sigma$ determines the
width of the Gaussian envelop and often varies in association
with the grid spacing, i.e., $\sigma = r h$, where $r$ is a
parameter depends on the stencil parameter $M$. For given 
$M$, there is wide range of near optimal $r$ values to be selected. 
In the DSC algorithm, for derivative approximation, one first analytically 
differentiates the kernel
\begin{eqnarray}
\delta^{(1)}(x-x_k) & = & \Big\{ \frac{\cos 
\frac{\pi}{h}(x-x_{k})}{(x-x_{k})} - \frac{\sin
\frac{\pi}{h}(x-x_{k})}{\frac{\pi}{h}(x-x_{k})^{2}}
\nonumber\\ &  &
- \frac{\sin \frac{\pi}{h}(x-x_{k})}{\frac{\pi}
{h}\sigma^{2}} \Big\}
\exp \Big( -\frac{(x-x_{k})^{2}}{2\sigma^{2}} \Big),
\end{eqnarray}
\begin{eqnarray}
\delta^{(2)}(x-x_k) & = &
\Big\{ -\frac{\frac{\pi}{h} \sin \frac{\pi}{h}
(x-x_{k})}{(x-x_{k})} -2 \frac{\cos \frac{\pi}{h}(x-x_{k})}
{(x-x_{k})^{2}} - 2\frac{\cos \frac{\pi}{h}(x-x_{k})}
{\sigma^{2}} \nonumber\\ &  & 
+2\frac{\sin \frac{\pi}{h}(x-x_{k})}{\frac{\pi}{h}
(x-x_{k})^{3}}
+ \frac{\sin \frac{\pi}{h}(x-x_{k})} {\frac{\pi}{h}(x-x_{k})
\sigma^{2}} \nonumber\\ &  & 
+\frac{\sin \frac{\pi}{h}(x-x_{k})}{\frac{\pi}{h}\sigma^{4}}
(x-x_{k}) \Big\}
\exp \Big( -\frac{(x-x_{k})^{2}}{2\sigma^{2}} \Big).
\end{eqnarray}
If $x=x_{k}$, one takes the limit
\begin{equation}
\delta^{(1)}(0) =0, \quad \textrm{and} \quad
\delta^{(2)}(0) = -\frac{1}{3}
\frac{3 + \frac{\pi^{2}}{h^{2}} \sigma^{2}}{\sigma^{2}}.
\end{equation}
Then we generate the difference coefficients of Eq. (\ref{center}) by evaluating at
the point of differentiation $ \delta^{(n)}_j=\delta^{(n)}(x-x_j)|_x$.
It is noted that only when one differentiates at $x=0$, the DSC first derivative
coefficients may be rewritten in terms of the sum-acceleration difference
by taking
\begin{equation}
w_{Mk}=\exp( - a^2 k^2),
\end{equation}
where $a^2=\frac{1}{2 r^2}$.
However, the second order and  higher order DSC-RSK method derivative coefficients 
cannot be cast into such a form.

\end{itemize}

\section{Results and discussions}\label{results}

To examine ``Proof'''s   claims in detail, we first present a comparison of the FD, Boyd's 
FD, Euler, MEuler, Sech and DSC-RSK methods for differentiating $e^{(ikx)}$ in Section \ref{deriva}.
It is noted that  $e^{(ikx)}$ involves only a single wavenumber, and might not be practical in 
applications.  However, the subject was chosen by ``Proof''. In fact, most 
problems studied in this work involve a wide range of wavenumbers.
As discussed earlier, the performance of a 
numerical method in differentiating $e^{(ikx)}$ does not always directly translate
into its performance in differentiating other functions, and has less to do with 
its performance in solving a differential equation. Since the ultimate test
of numerical methods is their behavior in solving differential equations, it is 
therefore important to examine the above mentioned methods for 
boundary value,  initial value, eigenvalue and nonlinear problems. 
It is noted that a boundary value problem studied here  was in fact previously 
designed by Boyd for comparing the performances of the FD, Euler, and Sech  methods \cite{Boyd94}. 

Both the Sech and DSC-RSK  methods have an adjustable parameter. For the Sech method,
the parameter $D$ depends on the spectral property of the solution.
Thus, for a fixed problem, it is independent of the half band-width $M$,
while the parameter $r$ of the RSK depends on $M$. In the following 
studies, the numerically-tested optimal $r$ and $D$ will be chosen,
and will be reported. In boundary value problems,
exact values outside the computational domain are used to impose boundary conditions.

\subsection{Differentiating $e^{ikx}$} \label{deriva}

\begin{figure*}[!tb]
\begin{center}
(a) \psfig{figure=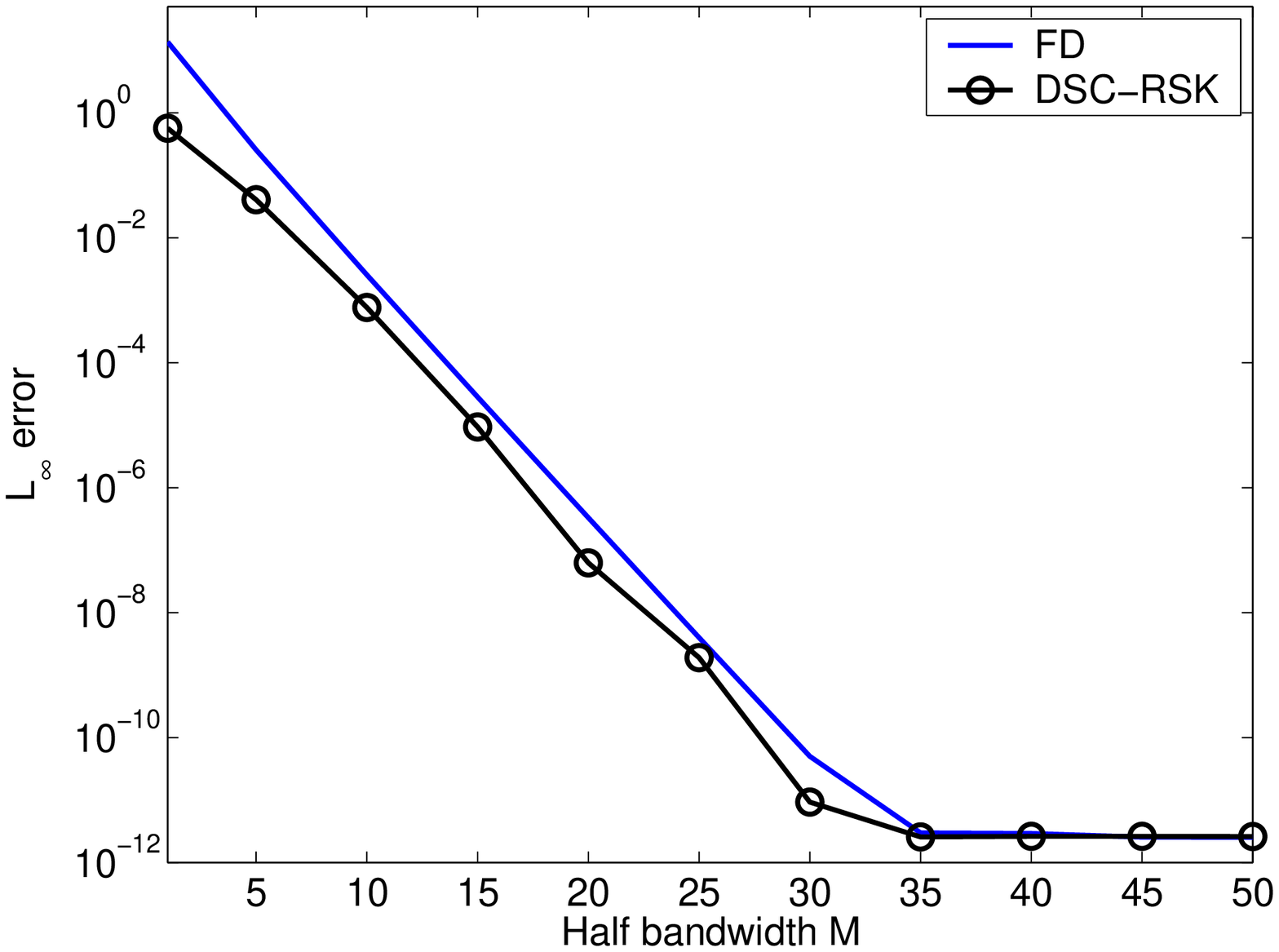,width=0.45\linewidth,height=6.5cm} 
(b) \psfig{figure=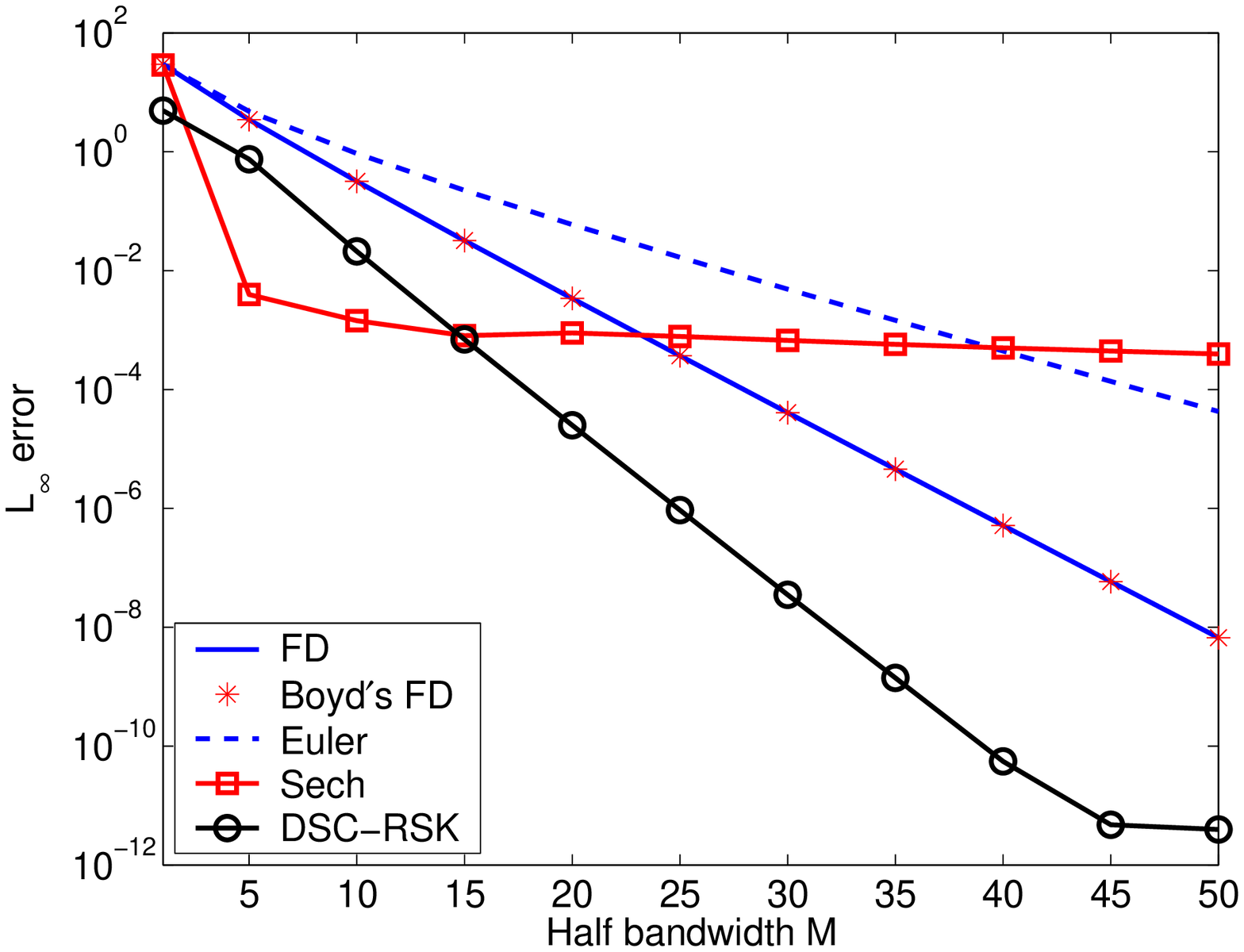,width=0.45\linewidth,height=6.5cm} \\
\end{center}
\caption{Error analysis in differentiating $\exp(ikx)$. 
(a) $kh=\frac{9}{20}\pi$;  (b) $kh=\frac{3}{5}\pi$.}
\label{fig.dtest} 
\end{figure*}

We first consider the approximation of differentiating $u(x)=\exp(ikx)$. 
It is claimed in ``Proof''  that the DSC is more accurate than the FD  only when
$kh>\frac{\pi}{2}$. A counterexample can be easily constructed at $kh=\frac{9}{20}\pi $
by considering a grid of 201 evenly spaced points in the interval of $x\in [0,2 \pi]$
and setting $k=45$. A near optimal DSC-RSK parameter $r$ is used in our computation for 
each $M$: $(M,r)=(1,1.2)$, (5,1.9), (10,2.7), (15,3.2),
(20,3.6), (25,4.0), (30,4.3), (35,4.7), (40,5.1), (45,5.5), and (50,5.9).
$L_{\infty}$ errors are depicted in Fig. \ref{fig.dtest} (a).
The numerical results clearly indicate that the DSC-RSK errors
are smaller than those of the FD method for a wavenumber less than $\frac{\pi}{2h}$, albeit 
the FD method is usually more accurate when $k$ is much smaller. This gives rise to our 
first counterexample of ``Proof''. 

Ironically, these results were partially implied in  Fig. 2 of ``Proof'', plotted with $M=33$. 
However, the results presented here indicate that the DSC-RSK method outperforms the finite difference 
for many small $M$ values, ranging from 1 to 35.

We next consider a problem with a medium-high wavenumber,
$kh=\frac{3}{5} \pi$, by increasing $k$ to 60 in the last case. 
This wavenumber is still far away from the aliasing limit. 
A near optimal Sech parameter, $D=0.17$ is used for all $M$ from 1 to 50,
while $r$ for the DSC-RSK method is chosen as: $(M,r)=(1,2.3)$, (5,2.4), (10,3.1), (15,3.7),
(20,4.2), (25,4.7), (30,5.1), (35,5.5), (40,5.8), (45,6.2), and (50,6.5).
$L_{\infty}$ errors are depicted in Fig. \ref{fig.dtest} (b).
The errors of the FD and Boyd's FD methods are identical.
The Sech method has a better performance only for $1<M<15$. Its accuracy 
does not improve as $M$ is increased beyond $M=15$. Therefore, if  $L_\infty$
error is required to be less than $10^{-4}$ for this problem, the Sech method cannot 
make it, whereas all other methods have a chance. The DSC-RSK method 
clearly outperforms the FD method up to 1000 times, and Boyd's spectrally-weighted 
finite difference (Sech)  method up to $10^8$ times at some large stencils.
This gives to another counterexample to ``Proof''.

\subsection{Differentiating a function having exponentially decaying amplitudes in wavenumbers $k$}\label{expon}

\begin{figure*}[!tb]
\begin{center}
(a) \psfig{figure=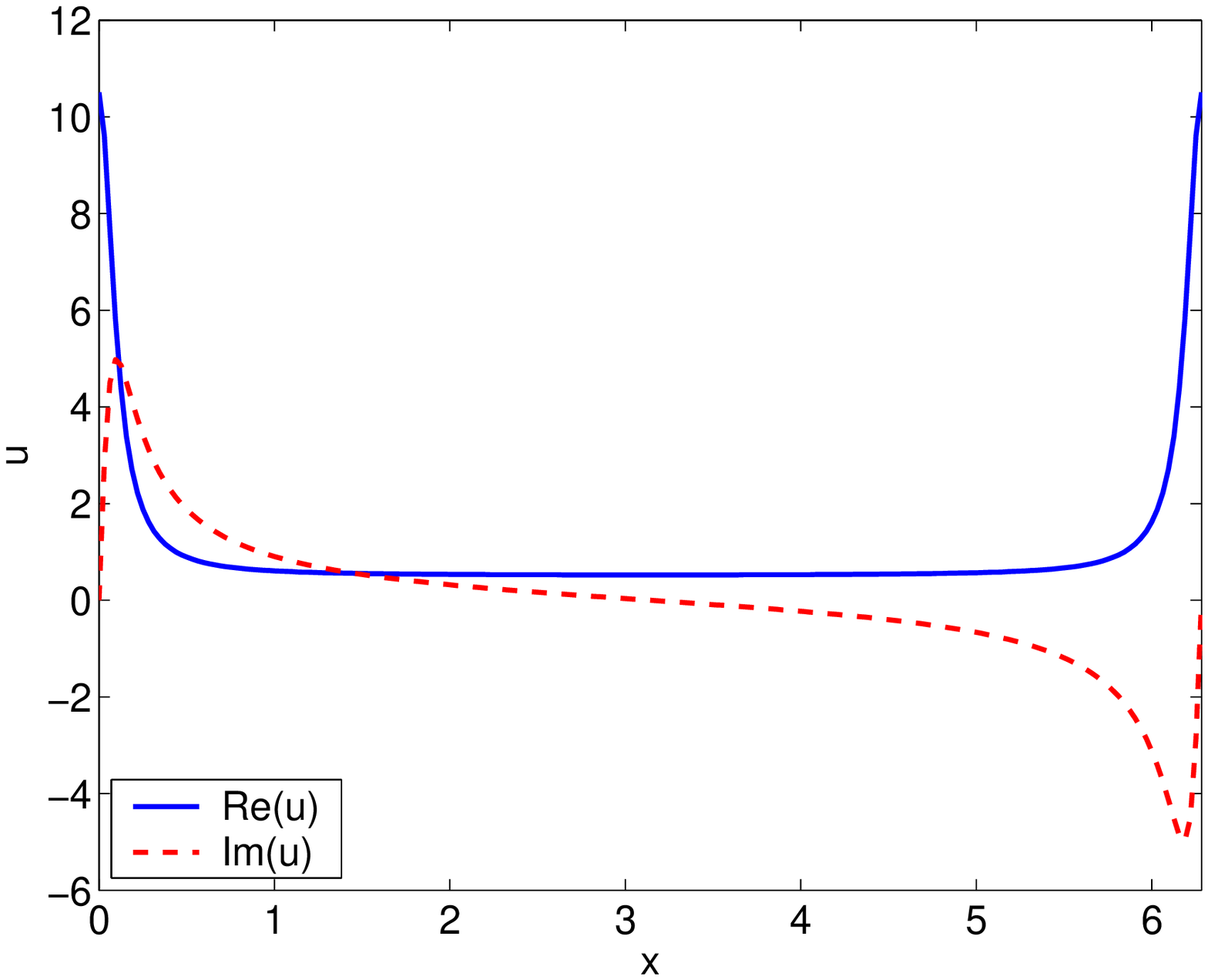,width=0.45\linewidth,height=6.5cm} 
(b) \psfig{figure=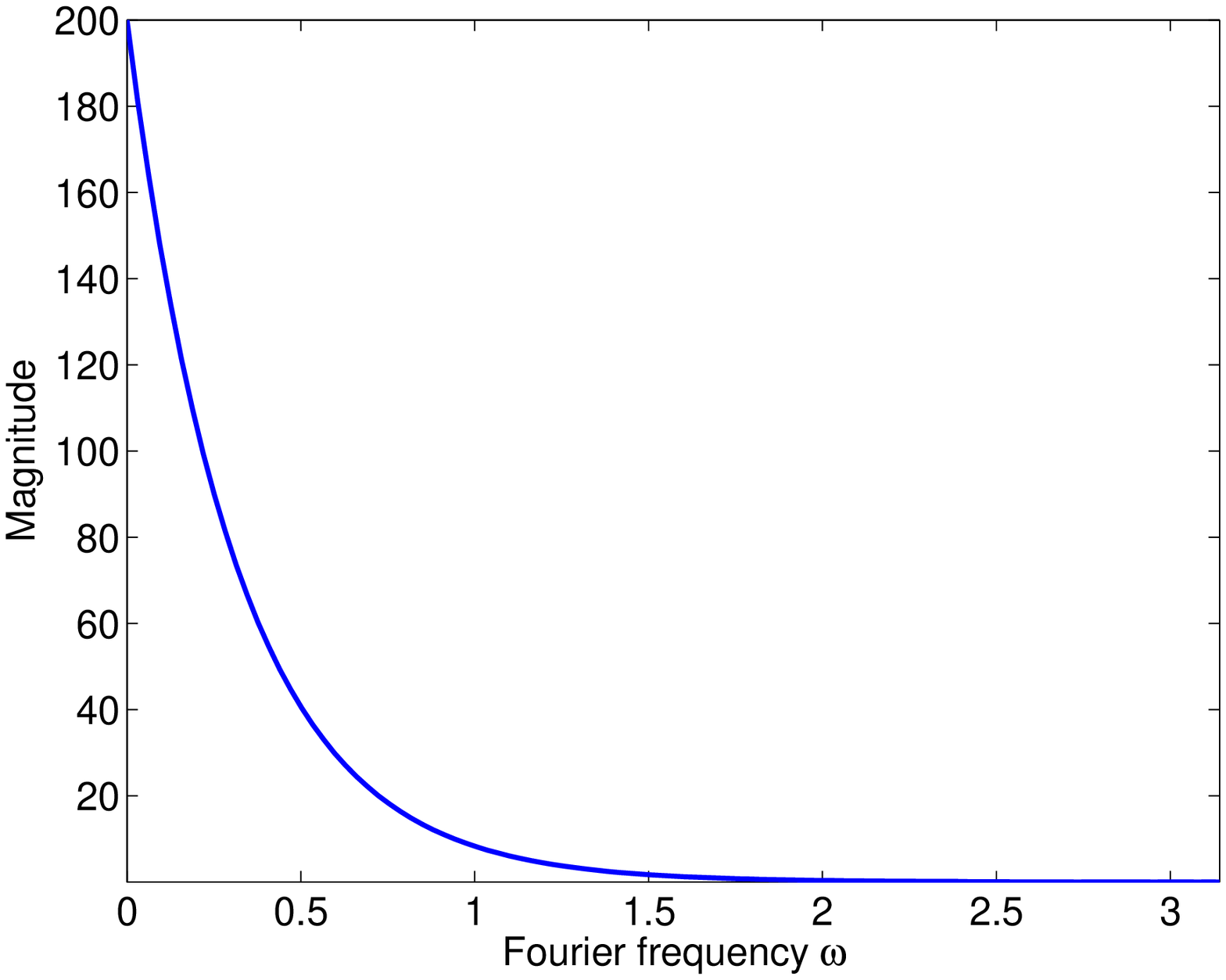,width=0.45\linewidth,height=6.5cm} \\
\end{center}
\caption{
(a) The function with exponentially decaying amplitudes in wavenumbers;
(b) Its Fourier frequency response.}
\label{fig.ftest1} 
\end{figure*}

\begin{figure*}[!tb]
\begin{center}
\psfig{figure=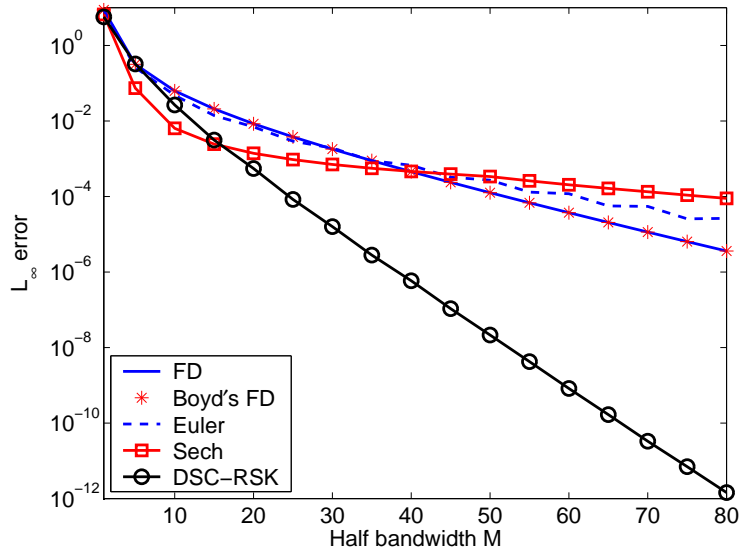,width=0.6\linewidth}  
\end{center}
\caption{Error analysis in differentiating the function having exponentially decaying amplitudes in 
wavenumbers.
}
\label{fig.ftest2} 
\end{figure*}

Although high frequency problems are common in science and engineering, some 
physical problems have their spectral distribution centered in the low frequency 
part, i.e., their amplitudes in large wavenumbers decay exponentially.  Therefore,
we consider the approximation of differentiating such a function, which
can be constructed as a weighted summation of Fourier basis functions $\exp(ikx)$
\begin{displaymath}
u(x)=\sum_{k=0}^{k_{\rm max}} \exp(-k \sigma) \exp( i k x), \quad x \in [0,2 \pi].
\end{displaymath}
Here, we choose $\sigma=0.1$ and $k_{\rm max}=80$. We carry out the differentiation on  a grid of 201 evenly 
spaced points. A near optimal Sech parameter, $D=0.17$ is used for all $M$ from 1 to 80,
while $r$ for the DSC-RSK method is chosen as: $(M,r)=(1,0.9)$, (5,1.8), (10,2.9), (15,3.8),
(20,4.6), (25,5.3), (30,5.9), (35,6.5), (40,7.1), (45,7.6), (50,8.1), (55,8.1),
(60,9.0), (65,9.4), (70,9.8), (75,10.3), and (80,10.6).

``Proof'' claims that ``Except for the very unusual case of low-pass filtered functions, 
that is, functions with negligible amplitude in small wavenumbers $k$, the DSC/LDAF is 
less accurate than the FD for all stencil widths $M$''. In the present problem, 
amplitudes in wavenumbers decay exponentially.  Fig. \ref{fig.ftest1}
gives the function and its Fourier frequency response. Fig. \ref{fig.ftest2} depicts
the error analysis of various numerical methods. Contrary to ``Proof'''s claim,  the DSC-RSK 
 method is more accurate than finite differences for all stencil widths $M$ 
considered. The Sech method performs better only when  $1<M<15$, and becomes the worst for 
large $M$ values. The DSC-RSK method clearly outperforms the FD method up to $10^6$ times, and the Sech 
 method up to $10^8$ times at some large stencils.

\subsection{ Boyd's boundary value problem }\label{boyds}

\begin{figure*}[!tb]
\begin{center}
\psfig{figure=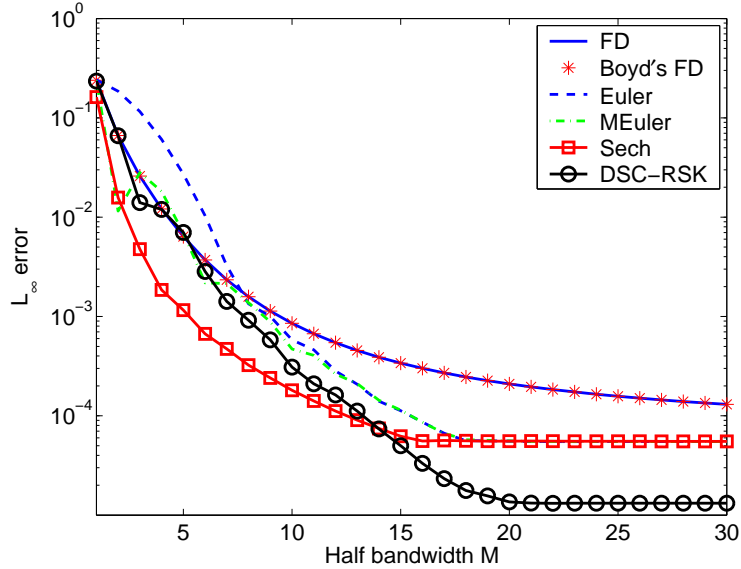,width=0.6\linewidth}  
\end{center}
\caption{Error analysis of Boyd's boundary value problem.}
\label{fig.bvp1} 
\end{figure*}
Boyd  employed a boundary value problem \cite{Boyd94} to demonstrate the superior 
performance of his spectrally-weighted difference, the Sech method,  over the FD scheme
\begin{eqnarray}
u_{xx} - u & = & f,  \nonumber \\
f & = &  \left\{ -2 \mbox{sech}^3(x)-\frac{\pi^2}{4 h^2}\mbox{sech} (x) \right\}
\cos (\frac{\pi x}{2 h})+\frac{\pi}{h}\mbox{sech} (x)\mbox{tanh} (x)
\sin (\frac{\pi x}{2 h}), \nonumber \\
u & = & \mbox{sech} (x) \cos (\frac{\pi x}{2 h}), 
\end{eqnarray}
where $h$ is the spacing. A grid of 201 evenly spaced points is used spanning the 
interval  $x\in [-30,30]$, which gives $h=0.3$. The  analytical solution outside
the computational domain was employed to impose boundary conditions.
Boyd examined the behavior of numerical methods by varying the stencil width $2M+1$ 
from $M=1$ to 20, while the results for the Sech method was only given from $M=1$ to 8. 
A near optimal Sech parameter, $D=0.25$ was used for all $M$ (see \cite{Boyd94}).
In this work, we shall solve the problem   
also extending 
the upper limit of the stencil parameter from $M=20$ to 30. A near optimal DSC-RSK parameter 
$r$ is used in our computations for each $M$: $(M,r)=(1,2.3)$, (2,1.4), (3,1.7), (4,2.2), 
(5,2.4), (6,2.5), (7,2.6), (8,2.9), (9,3.1), (10,3.3), (11,3.5), (12,3.7), (13,4.0), 
(14,4.1), (15,4.3), (16,4.5), (17,4.7), (18,4.9), (19,5.1), 
(20,5.4), (21,5.4), (22,5.4), (23,5.4), (24,5.4), (25,5.4), 
(26,5.4), (27,5.4), (28,5.4), (29,5.4), and (30,5.4).

\begin{figure*}[!tb] 
\begin{center}
(a) \psfig{figure=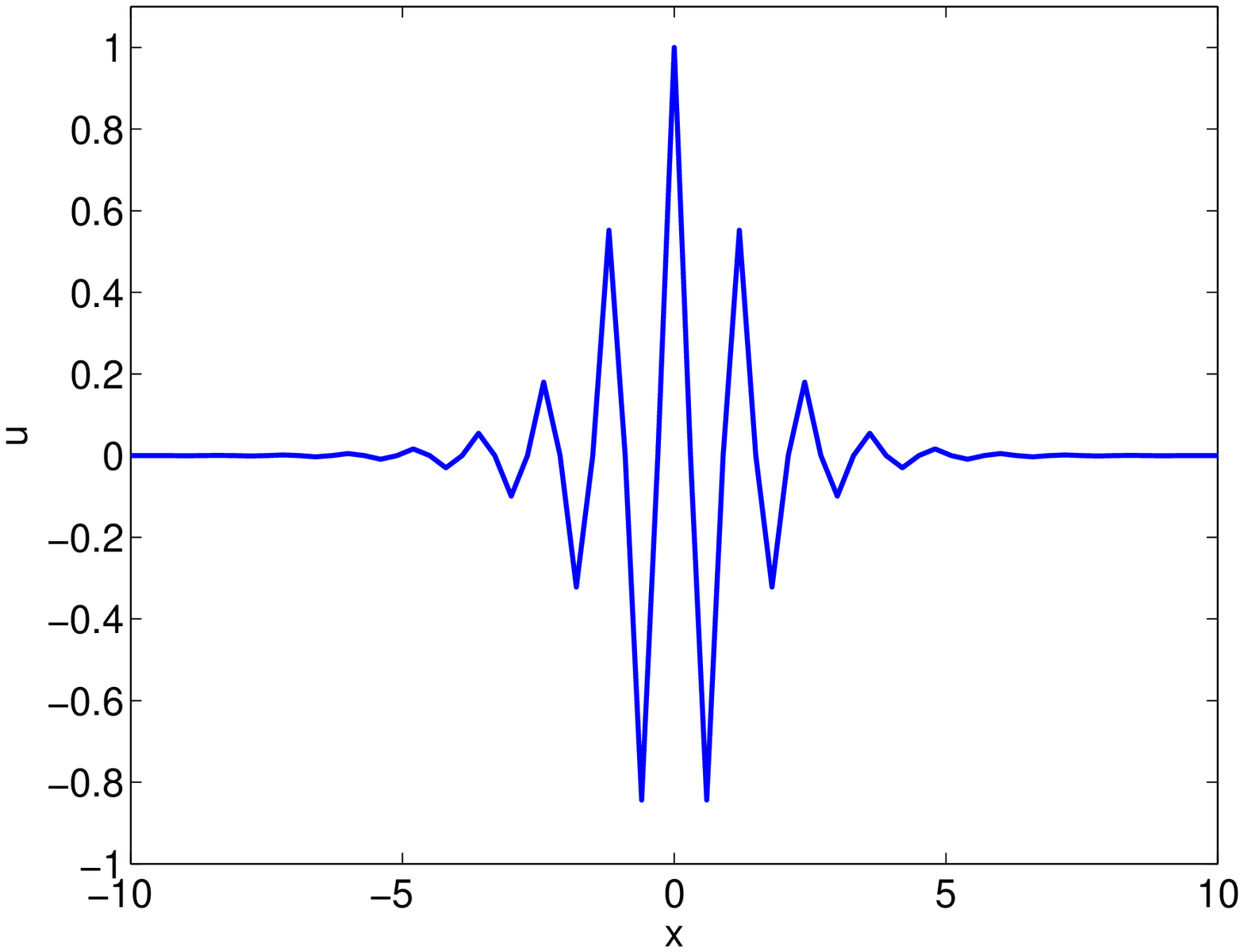,width=0.45\linewidth,height=6.5cm} 
(b) \psfig{figure=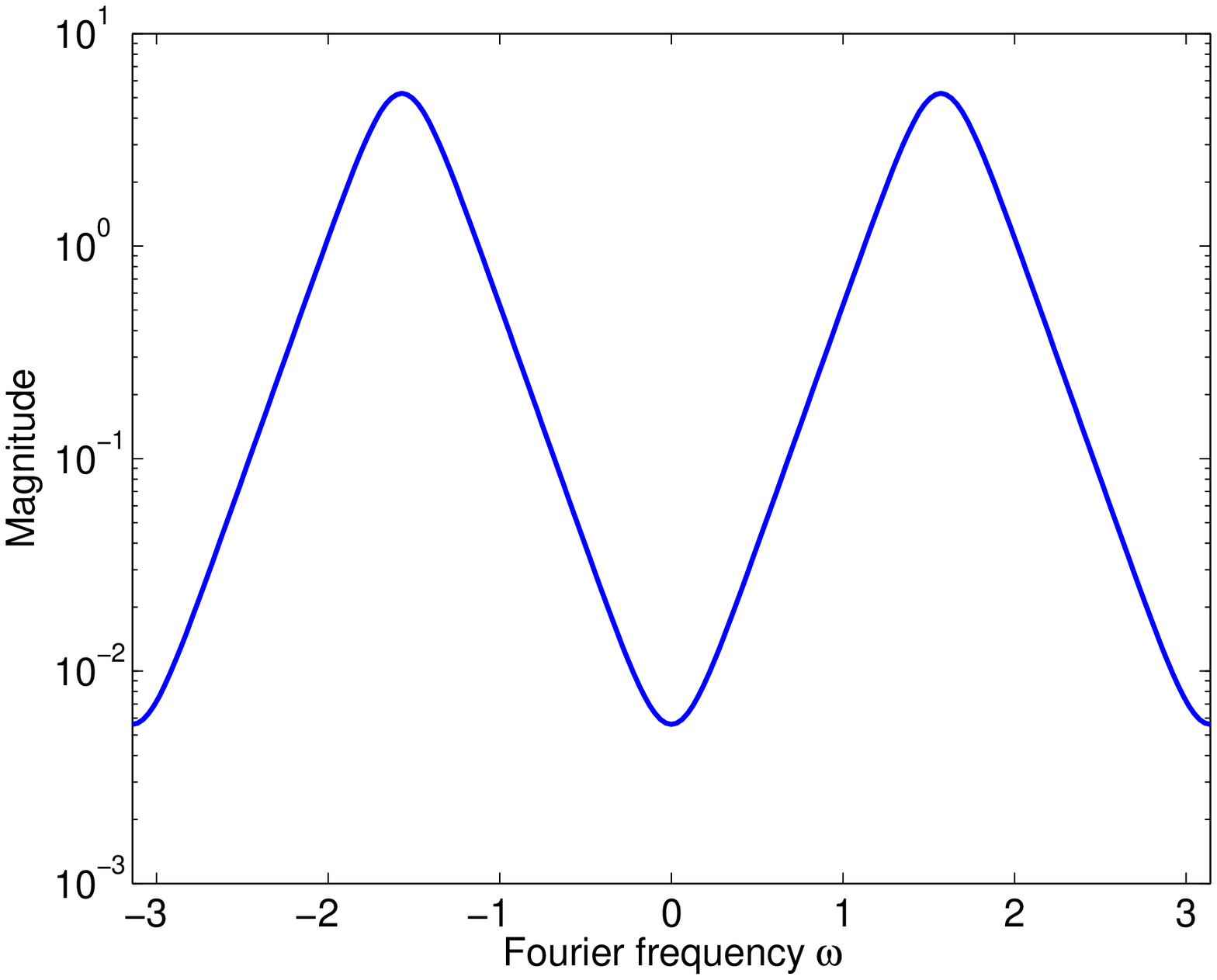,width=0.45\linewidth,height=6.5cm} \\
\end{center}
\caption{(a) Solution of Boyd's boundary value problem; 
         (b) Its Fourier frequency response.
}
\label{fig.BVPsolu} 
\end{figure*}

$L_{\infty}$ errors are depicted in Fig. \ref{fig.bvp1}. We note that the behaviors of FD, 
Euler, and Sech  methods are identical to those in  Fig. 4 of Ref. \cite{Boyd94}. 
Indeed, for $M<13$ the Sech method performs better than the FD and DSC-RSK schemes. However, for 
a large stencil width, the DSC-RSK method clearly outperforms all other methods, including the FD 
and Sech  methods.  Therefore, Boyd's boundary value problem is a counterexample of his claims in ``Proof''.

The flat and comparatively large errors  of  all six different methods for $M\geq20$ 
are suspicious, indicating that the errors are limited by the original design of the problem. 
To investigate further, we analyzed the Fourier frequency response of the exact solution
as shown in Fig. \ref{fig.BVPsolu} (b). Indeed, the magnitude of the Fourier frequency response 
has a large truncation at $|\pi|$, which gives rise to inherent aliasing errors.  Therefore, this 
problem is not appropriate for demonstrating the full potential of different numerical methods.

\subsection{Boundary value problem with a confined  distribution of wavenumbers}\label{bound1}

\begin{figure*}[!tb]
\begin{center}
(a) \psfig{figure=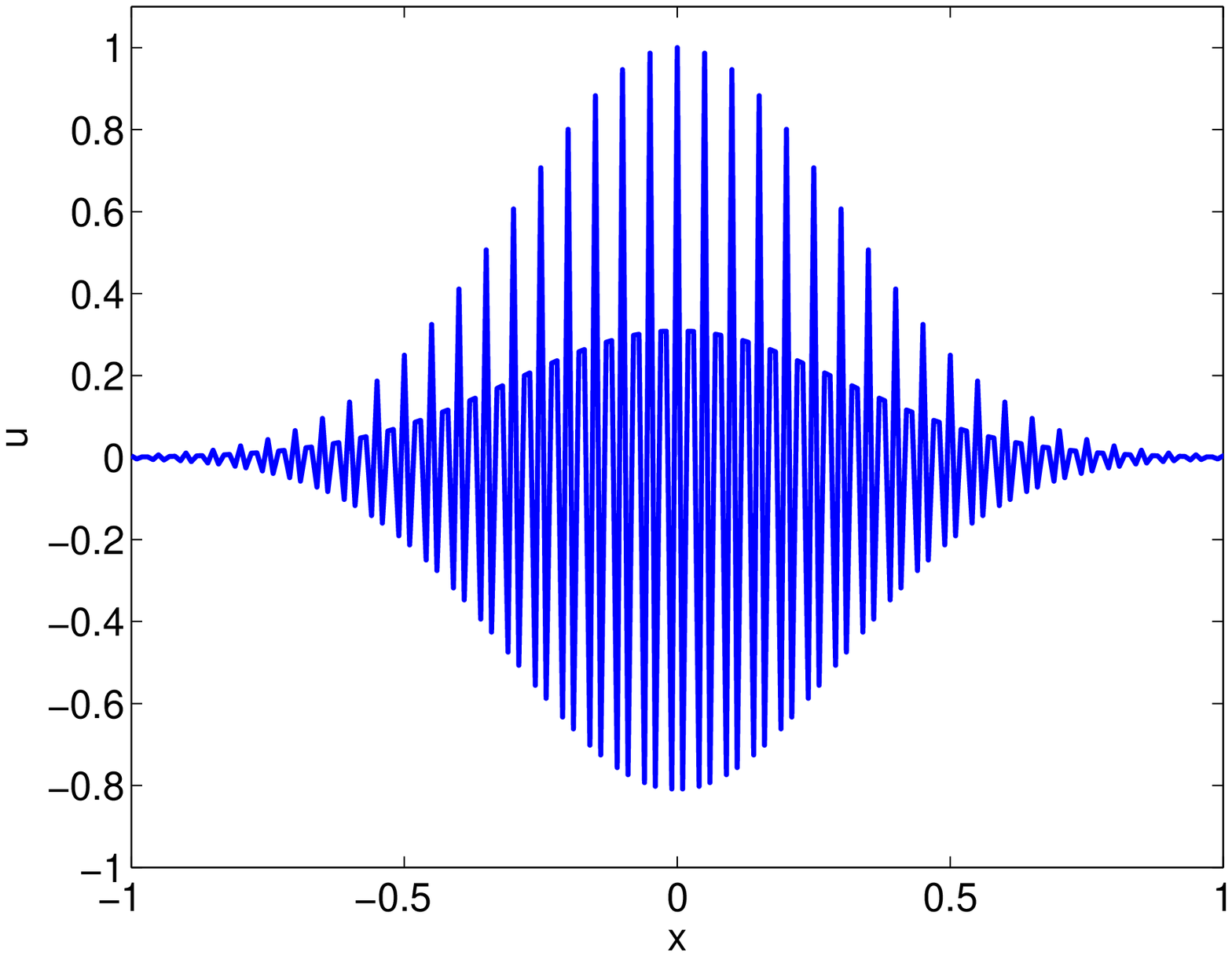,width=0.45\linewidth,height=6.5cm} 
(b) \psfig{figure=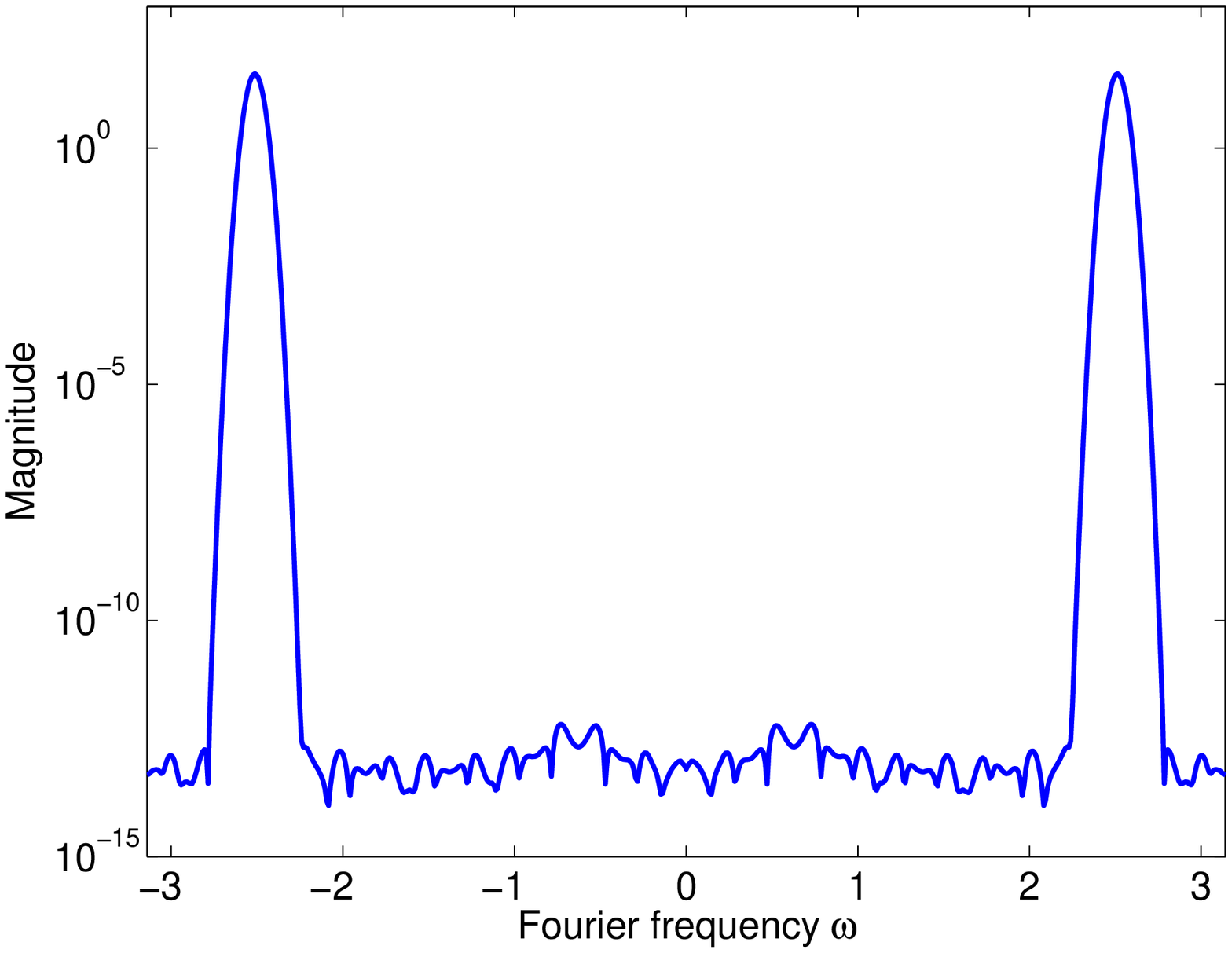,width=0.45\linewidth,height=6.5cm} \\
\end{center}
\caption{ 
(a) The solution of a boundary value problem with a confined  distribution of wavenumbers;
(b) Its Fourier frequency response.}
\label{fig.BVPsolu2} 
\end{figure*}

Since the above boundary value problem, taken from \cite{Boyd94},
 cannot reveal the true  accuracy 
of numerical methods, let us consider a modification 
\begin{eqnarray}
u_{xx} - u & = & f, \quad x \in [-3,3], \nonumber \\
u & = & \exp(-\frac{x^2}{2 a^2}) \cos(b \pi x), \\
f & = &  \frac{2 b \pi x}{a^2}\exp(-\frac{x^2}{2 a^2}) \sin(b \pi x)
+\left(\frac{x^2}{a^4}-\frac{1}{a^2}-(b\pi)^2-1 \right)
\exp(-\frac{x^2}{2 a^2}) \cos(b \pi x). \nonumber
\end{eqnarray}
By choosing $a=0.3$ and $b=80$, the solution is oscillatory and its frequency 
response shows a peak around $\pm 2.5$,  see Fig. \ref{fig.BVPsolu2}. 
Note that the frequency magnitude is as small as $10^{-14}$ near $\pm \pi$. Therefore,
in contrast to the previous  boundary value problem, the present problem has negligible 
inherent aliasing error. Due to the frequency distribution in the large wavenumber
region, this  case  favors numerical methods that perform well for large wavenumbers. 
 We employ a grid of $N=601$ points in the interval $x\in [-3,3]$. 
In the Sech method, the nearly optimized $D=0.18$ is used. Other $D$ could be slightly 
better but the amount of improvement would be not significant enough to alter the conclusion of our 
analysis. For the DSC-RSK method, near optimal $r$ values are chosen as 
$(M,r)=(1,2.7)$, (5,3.9),  (10,4.5), (15,6.0), (20,6.4), (25,7.1), (30,7.6), (35,8.1), (40,8.6),
(45,9.2), (50,9.6), (55,10.1), (60,10.5), (65,11.0), (70,11.3), (75,11.8), and (80,12.3).

\begin{figure*}[!tb]
\begin{center}
\psfig{figure=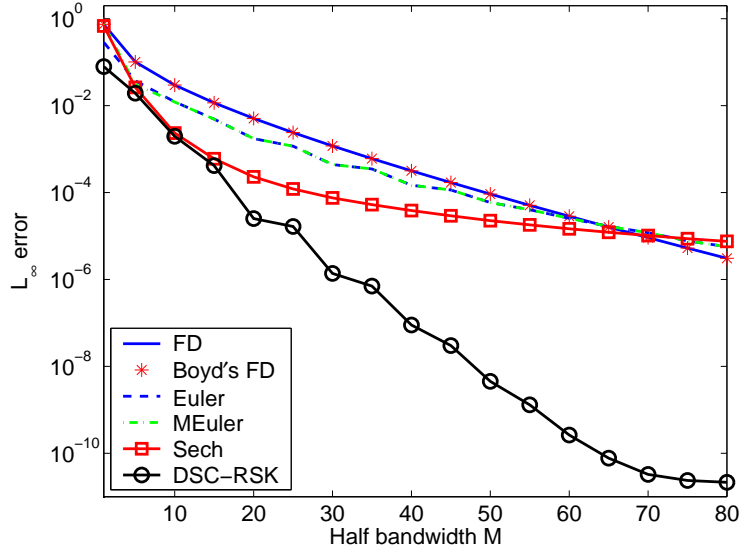,width=0.6\linewidth}  
\end{center}
\caption{Error analysis for the boundary value problem with confined wavenumber distribution.}
\label{fig.bvp2} 
\end{figure*}

Since this is a problem involving large  wavenumbers, according to ``Proof'', the 
spectrally-weighted differences are designed for handling this kind of problems, 
and will outperform the DSC-RSK method. 
Fig. \ref{fig.bvp2} illustrates the $L_{\infty}$ errors for six different 
numerical methods for this case. 
It is interesting to note that Boyd's FD method yields the same results as the standard 
FD method for all $M$ considered. The MEuler method is slightly more accurate than the Euler method when $M$ is relatively
small, while when $M > 10$, their results are identical.
The Sech method is more accurate than these four differencing methods over a wide range of $M$
for this high frequency problem. 
However, it is clear from Fig. \ref{fig.bvp2} that the DSC-RSK method is the most accurate method for all
$M$ values. 
In particular, when $M=80$, the DSC-RSK method is about five orders more accurate than any other numerical method
considered. 
In fact, at large $M$, the accuracy of the DSC-RSK method, $L_{\infty}=10^{-11}$, is limited 
by the iterative algebraic solver, the standard preconditioned biconjugate gradient (PBCG) 
method \cite{Press}.

\subsection{Boundary value problem with a wide range of wavenumbers }\label{bound2}

\begin{figure*}[!tb]
\begin{center}
(a) \psfig{figure=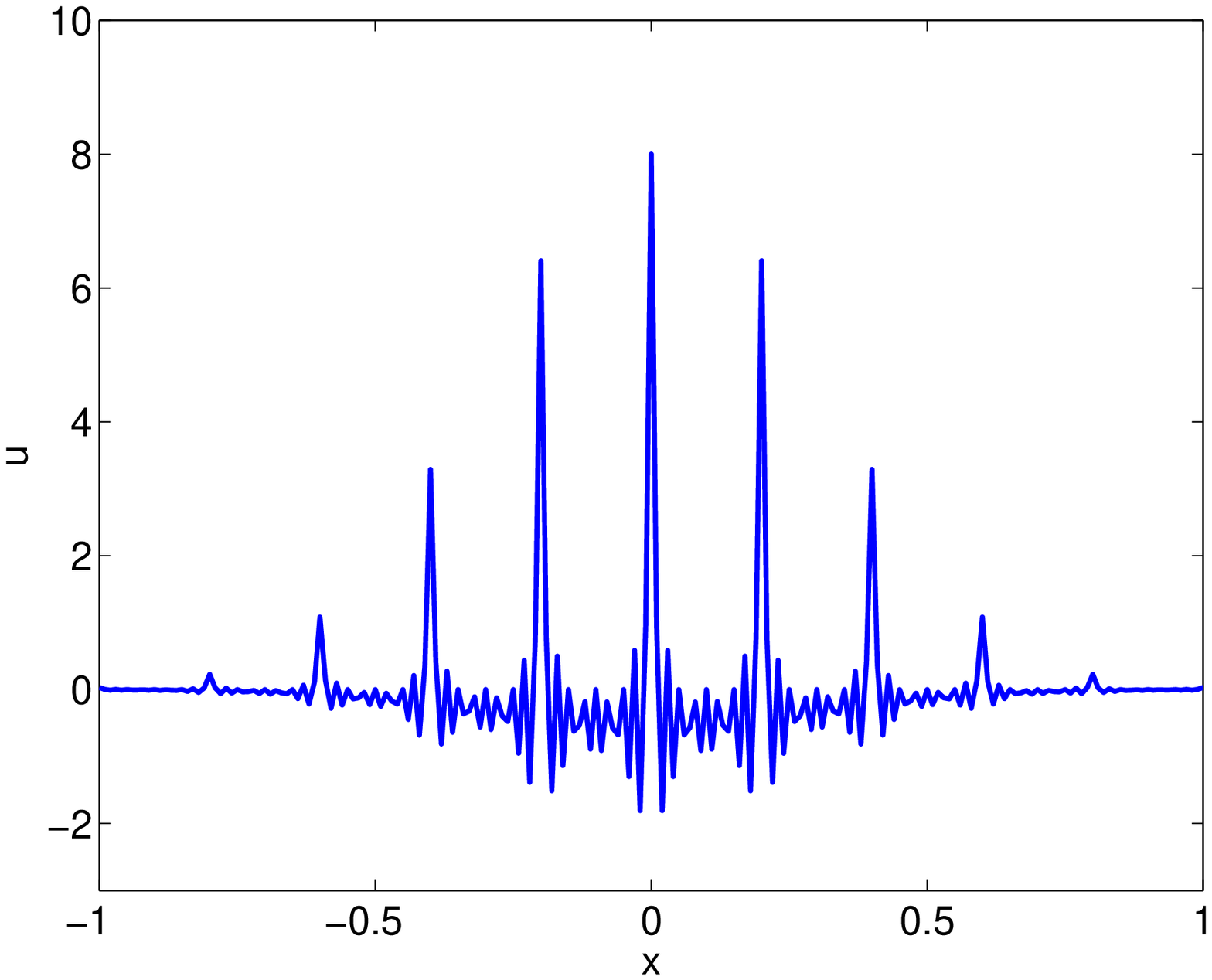,width=0.45\linewidth,height=6.5cm} 
(b) \psfig{figure=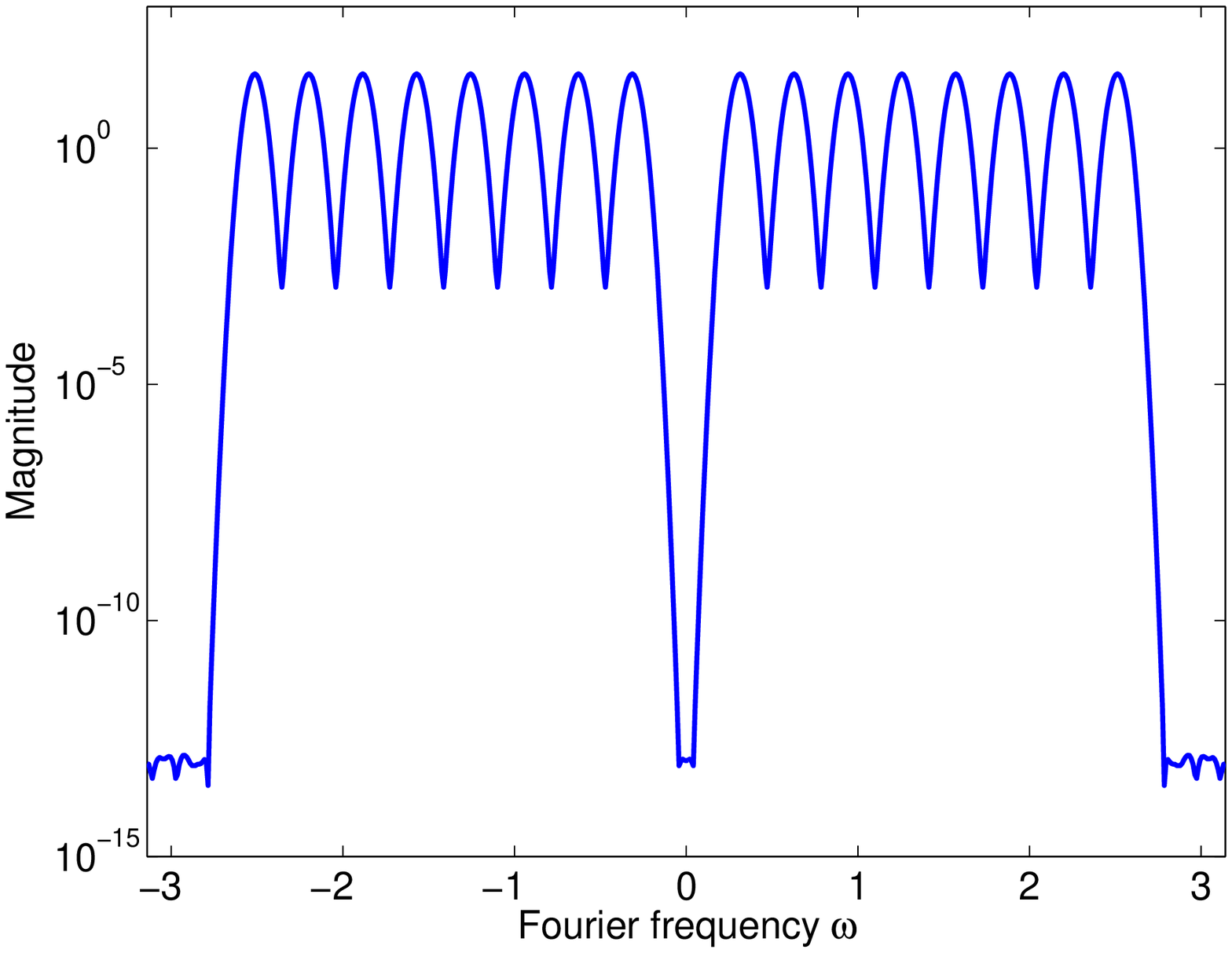,width=0.45\linewidth,height=6.5cm} \\
\end{center}
\caption{(a) The solution of a boundary value problem with a wide range of wavenumbers; 
(b) Its Fourier frequency response.
}
\label{fig.BVPsolu3} 
\end{figure*}

It is important to study problems with a wide range of wavenumbers.
Let us  consider a boundary value problem
\begin{eqnarray}
u_{xx} - u & = & f, \quad x \in [-3,3], \nonumber \\
u & = & \sum_{b=1}^{8} \exp(-\frac{x^2}{2 a^2}) \cos(10 b \pi x), \\
f & = & \sum_{b=1}^{8}  \frac{20 b \pi x}{a^2}
\exp(-\frac{x^2}{2 a^2}) \sin(10b \pi x)
+\sum_{b=1}^{8}  \left(\frac{x^2}{a^4}-\frac{1}{a^2}-(10b\pi)^2-1 \right)
\exp(-\frac{x^2}{2 a^2}) \cos(10b \pi x), \nonumber
\end{eqnarray}
where $a=0.3$.  
The solution and its frequency response are depicted in 
Fig.  \ref{fig.BVPsolu3}. It is seen that the frequency response has a wide range of 
wavenumbers. We have arranged the amplitudes in small wavenumbers to be as 
large as those in large wavenumbers (see Figs.  \ref{fig.BVPsolu3} (b)).

We still solve the problem with $N=601$. For the Sech method, it is not 
easy to find the optimal parameter $D$ for this problem due to a wide range 
of involved wavenumbers.  A near optimized 
value, $D=0.18$, is employed.  For the DSC-RSK method, near optimal parameters
are as follows:
$(M,r)=(1,1.1)$, (5,2.6),  (10,3.1), (15,4.6), (20,5.3), (25,6.1), (30,6.8), (35,7.5), (40,8.1),
(45,8.5), (50,9.0), (55,9.5), (60,10.0), (65,10.4), (70,10.8), (75,11.4), and (80,11.8).

\begin{figure*}[!tb]
\begin{center}
\psfig{figure=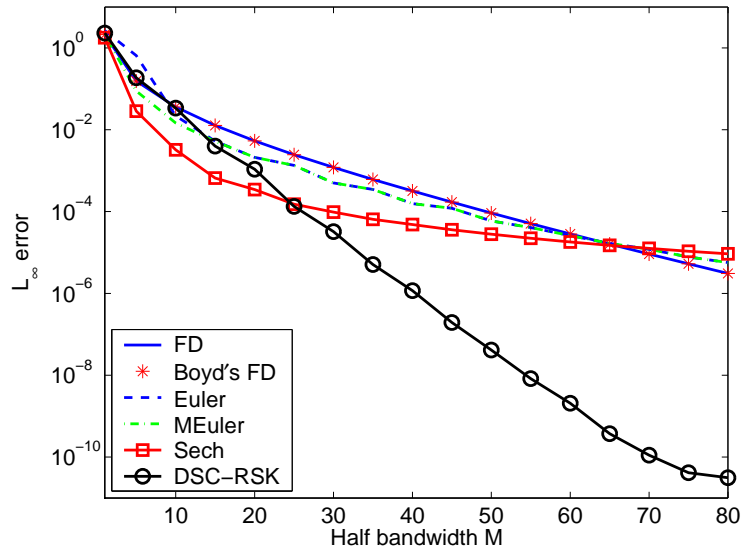,width=0.6\linewidth}  
\end{center}
\caption{Error analysis for the boundary value problem with a wide range of wavenumbers
}
\label{fig.bvp3} 
\end{figure*}

The $L_{\infty}$ errors are given in Fig. \ref{fig.bvp3} for different $M$.
The MEuler method is more accurate than the Euler method for small $M$, while when 
$M > 15$, their results are identical.
For a small $M$, such as up to $80$ in this example, Boyd's
FD method yields the same results as the FD method. However, when
$M>80$, Boyd's FD method is less accurate, probably because for large  $M$,
calculating  the weights of Boyd's FD method produces numerical errors.
The results of the Sech method  are similar to those in the last case.
In particular, the Sech method is more accurate than the FD, Boyd's FD, Euler and 
MEuler methods when $M < 70$, while 
when $M$ is sufficient large, the accuracy of the Sech method  is worse than the others. 
It is noted that the Sech method performs slightly better than the 
DSC-RSK method from $M < 25$, due to large amount of the small wavenumber components.
However, for a large $M$, the DSC-RSK method is much more accurate.

``Proof''  claims that ``For $a < a_{\rm  FD}$, the DSC errors 
are less than finite differences for $k$ near the aliasing limit, but much, much worse for smaller $k$. 
Except for the very unusual case of low-pass filtered functions, that is, functions with negligible 
amplitude in small wavenumbers $k$, the DSC/LDAF is less accurate than finite differences for all 
stencil widths $M$''.  
The present example satisfies all of  ``Proof'''s descriptions, except for the 
conclusion: the DSC-RSK method is more accurate than the FD method for all stencil widths examined.
When $M=80$, the DSC-RSK method outperforms the FD scheme and all other methods by more than
100,000 times.
If  high accuracy is desirable for this problem, the DSC-RSK scheme is the method of choice
among the six methods.

\subsection{Helmholtz equation with a constant source term}\label{const}

We next consider the Helmholtz equation with high wavenumbers $k$
\begin{eqnarray}
-u_{xx}(x)- k^2 u(x) & = & 1, \quad x \in [0,1]
\nonumber\\ 
\label{Helmholtz} u(0) & = & 0, \\
u_x(1)-i k u(1) & = & 0, \nonumber \\
u(x) & = & \frac{1}{k^2} \left[(1-\cos(kx)-\sin (k) \sin (kx)) + i
(\cos(k)-1)\sin(kx)\right]. \nonumber
\end{eqnarray}
This problem has its origin from electromagnetics and is notoriously 
challenging when  $k$ is large for a given number of grid points 
(see \cite{Babuska,Babuska2,Deraemaeker,Ihlenburg1,Ihlenburg2}). 
High wavenumbers are common in electrically-large systems.  
As ``Proof''  claims that spectrally-weighted differences are superior to the DSC for 
differentiating functions with large wavenumbers, we employ the Sech method.
We shall explore the capability of each high order 
method for resolving high frequency problems with a small mesh density close to 
the Nyquist limit.  A method that can achieve a spatial resolution near the Nyquist limit,
2 points per wavelength (PPW), without invoking much alias errors, is  more efficient 
for solving large scale problems. The DSC-RSK method was applied to this problem in our previous 
study \cite{BaWeZh04}. The objective of the present work is to examine ``Proof'''s   claim and 
to explore different numerical methods.

\begin{figure*}[!tb]
\begin{center}
(a) \psfig{figure=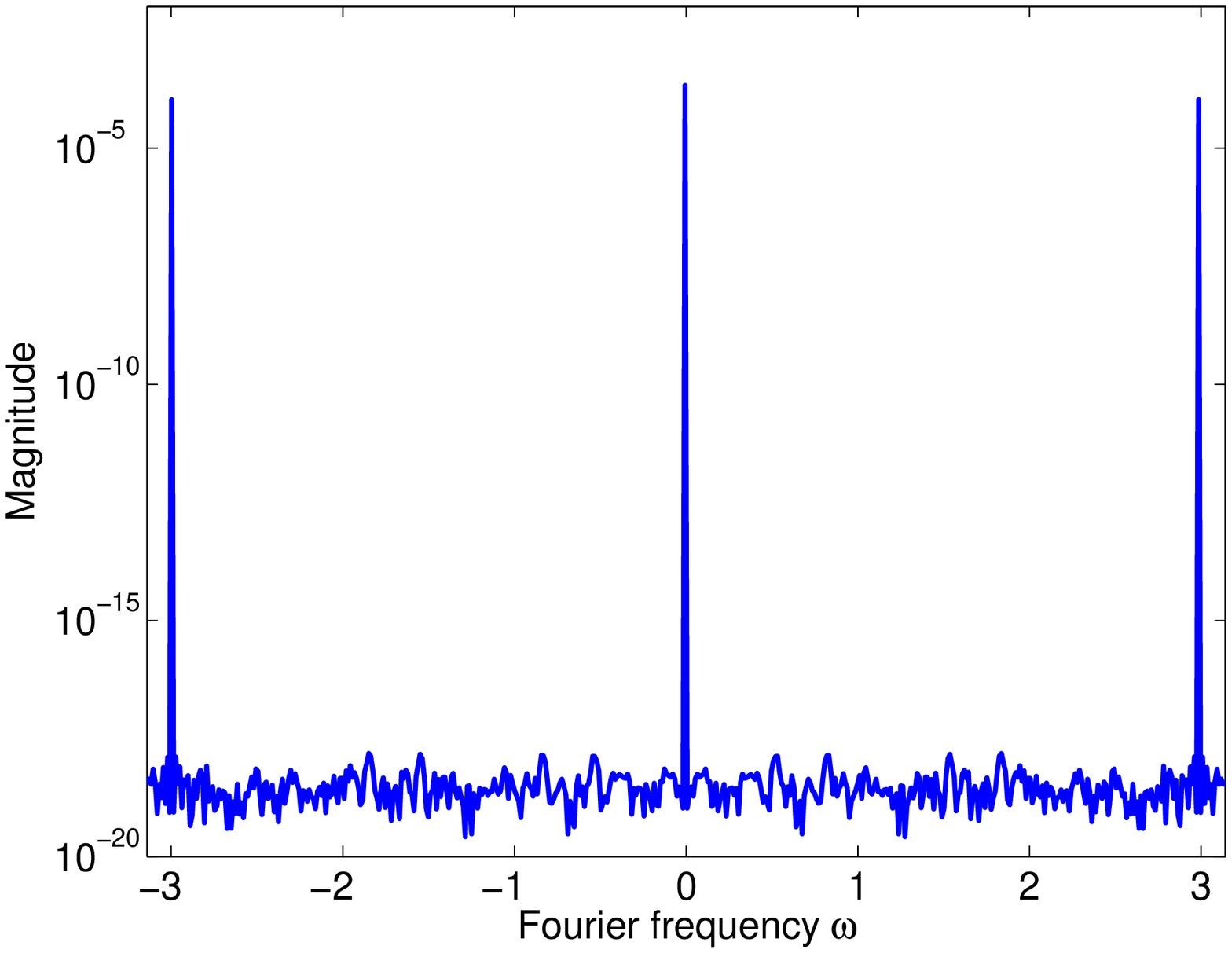,width=0.45\linewidth,height=6.5cm} 
(b) \psfig{figure=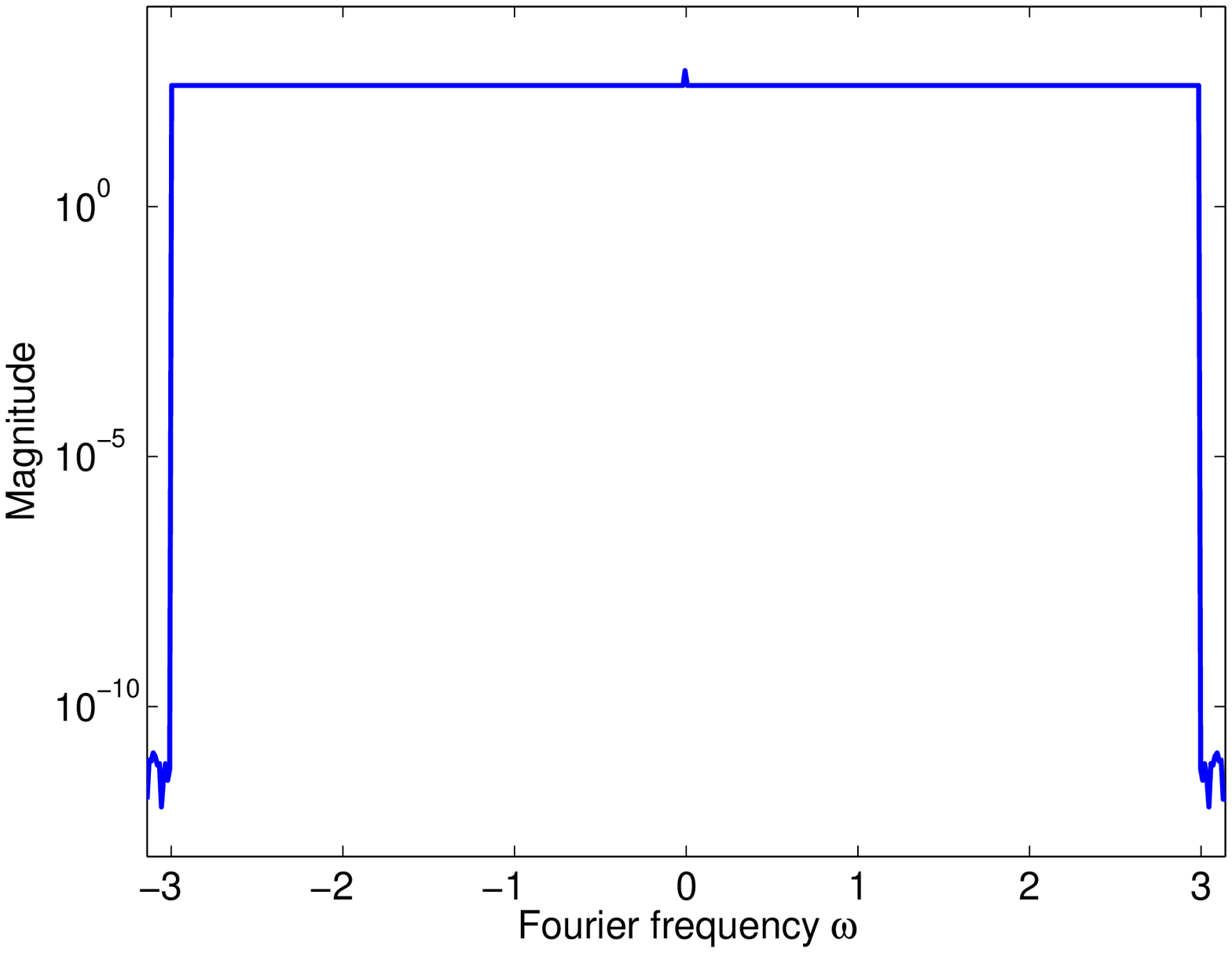,width=0.45\linewidth,height=6.5cm} \\
\end{center}
\caption{Fourier frequency response of the solution to the Helmholtz equation with high wavenumbers.
(a) With a constant source; 
(b) With a multiple frequency source.  
}
\label{fig.Helmsolu} 
\end{figure*}

We first consider the case of $k=500 \pi$. 
The Fourier frequency response of the solution is depicted  in Fig. \ref{fig.Helmsolu}.
We use a grid of $N=526$ points in the domain so that the grid density is about  
2.1 PPW, which is very close to the Nyquist sampling limit. Therefore,
the setting of problem is extremely difficult for numerical methods that do not work
well for high wavenumbers. The $H^1$-seminorm will be used to measure the error \cite{BaWeZh04}
\begin{equation}
e_1 :=\frac{|u-u_h|_1}{|u|_1},
\end{equation}
where $u_h$ is the numerical approximation to the solution $u$, and
$|u|_1=||\frac{d u}{d x}||$. 
Note that for this problem, large stencils are required to deliver reasonable results. 
We solve the problem with a variety of $M$ values ranging from  
$M=50$ to 500. Since  the values of $M$ are quite large, the total execution time for generating 
the first and second derivative weights is also reported for each method for a comparison.  
In the Sech method, a near optimized parameter $D=0.28$ is used for all $M$.
In the DSC-RSK method, we choose $(M,r)=(50,9.0)$, (100,27.9), (150,35.1), (200,38.1),
(250,42.3), (300,46.9), (350,49.6), (400,53.5), (450,56.0), and (500,60.0).

Since the numerical errors vary over 160 orders, it is inconvenient 
to put them in a figure. We therefore present them in Table \ref{table.Helm1}. 
It is seen from the table that for the present problem with PPW=2.1, the FD method 
does not improve much when $M$ is increased from 50 to 500. Its errors are 
intolerably large for all $M$ values. Therefore, the FD method is not efficient
for this highly oscillatory problem. 

\begin{table}[!tb]       
\caption{Numerical results and CPU times in seconds 
for solving the Helmholtz equation with a constant source term.
Here NaN and --- stand for ``not a number'' and ``not available'', respectively.
}
\label{table.Helm1}      
\begin{center}      
\begin{tabular}{lllllll}
\hline
\hline
 & \multicolumn{2}{c}{FD} & \multicolumn{2}{c}{Boyd's FD} & 
\multicolumn{2}{c}{Sech} \\
\hline
$M$ & Error & CPU  & Error & CPU  & Error & CPU \\
\hline
50 &   $9.98(-1)$ & $8.63(-4)$ & $9.98(-1)$ & $4.49(-4)$ & $5.71(-2)$ & 3.983\\
100 &  $2.32(+0)$ & $3.13(-3)$ & NaN & $3.41(-4)$        & $5.24(-1)$ & 15.69\\
150 &  $1.17(+0)$ & $7.02(-3)$ & --- & ---               & $8.88(-2)$ & 35.25\\
200 &  $1.36(+0)$ & $1.38(-2)$ & --- & ---               & $3.98(-2)$ & 62.77\\
250 &  $4.27(+0)$ & $2.03(-2)$ & --- & ---               & $4.85(-2)$ & 98.26\\
300 &  $6.57(-1)$ & $2.84(-2)$ & --- & ---               & ---  & --- \\
350 &  $5.09(-1)$ & $3.88(-2)$ & --- & ---               & ---  & --- \\
400 &  $4.10(-1)$ & $4.91(-2)$ & --- & ---               & ---  & --- \\
450 &  $3.23(-1)$ & $6.17(-2)$ & --- & ---               & ---  & --- \\
500 &  $2.49(-1)$ & $7.58(-2)$ & --- & ---               & ---  & --- \\
\hline
 & \multicolumn{2}{c}{Euler} & \multicolumn{2}{c}{MEuler} & 
\multicolumn{2}{c}{RSK}\\
\hline
$M$ & Error & CPU  & Error & CPU  & Error & CPU \\
\hline
50 &  $2.72(-1)$ & $4.77(-4)$ & $2.72(-1)$ & $4.72(-4)$  &  $1.16(+0)$ & $3.06(-4)$ \\
100 & $2.94(+1)$ & $4.01(-4)$ & $2.94(+1)$ & $4.08(-4)$  &  $1.12(-2)$ & $2.31(-4)$ \\
150 & $3.02(-1)$ & $5.95(-4)$ & $3.02(-1)$ & $6.16(-4)$  &  $7.52(-5)$ & $3.39(-4)$ \\
200 & $5.16(+58)$ & $5.83(-4)$ & $5.16(+58)$ & $7.30(-4)$  &  $2.27(-8)$ & $4.20(-4)$ \\
250 & $5.81(+73)$ & $6.86(-4)$ & $5.81(+73)$ & $7.02(-4)$  &  $7.34(-10)$ & $5.55(-4)$ \\
300 & $6.55(+88)$ & $2.67(-2)$ & $6.55(+88)$ & $2.67(-2)$  &  $7.08(-10)$ & $6.57(-4)$ \\
350 & $7.37(+103)$ & $2.04(-2)$ & $7.37(+103)$ & $2.08(-2)$  &  $3.45(-11)$ & $7.15(-4)$ \\
400 & $8.30(+118)$ & $2.56(-2)$ & $8.30(+118)$ & $2.65(-2)$  &  $2.83(-13)$ & $8.09(-4)$ \\
450 & $9.35(+133)$ & $2.74(-2)$ & $9.35(+133)$ & $2.82(-2)$  &  $1.50(-13)$ & $9.05(-4)$ \\
500 & $1.05(+149)$ & $3.64(-2)$ & $1.05(+149)$ & $3.64(-2)$  &  $1.77(-13)$ & $9.97(-4)$ \\
\hline
\hline
\end{tabular}       
\end{center}      
\end{table}

When $M = 50$,  Boyd's FD scheme yields the same results as the standard FD method. However, when $M = 100$, 
although the weights are generated in 0.0003 second, the rest of the computation lasts a 
very long time without delivering any reasonable result. It gets worse  for $M > 100$. 
This seems to suggest that Boyd's FD method is unstable for large $M$. 

Although the Sech method was proposed to deal with high wavenumbers, it does not work the way
as promised for very high wavenumber problems. One major problem with the 
Sech method is that it requires excessive CPU time for generating the weights. 
Since the Gram matrix is poorly conditioned, the SVD solver is 
often used, which is CPU demanding. Thus, the Sech method seems to be 
impractical for high frequency problems. 
Moreover, the SVD solver \cite{Press} fails to converge when $M \ge 300$ for the Sech method, 
probably because the matrix is seriously ill-conditioned for  large $M$. 
Except for $M=50$, the Sech method was less accurate than the DSC-RSK method for this challenging problem. 
At $M=250$, the DSC-RSK method is about $10^{8}$ times more accurate than the Sech method. If 
the minimal accuracy requirement is $10^{-3}$, no other methods except for the DSC-RSK method
can be used for this problem.  

The Euler and MEuler methods produce identical results. Their errors reach $10^{149}$
at $M=500$, indicating the instability in computing the factorial. Nevertheless,
the same sets of weights seem to work better in the next problem with a wide range of
high wavenumbers. Therefore, these two methods are not stable for being used in 
high wavenumber problems from the Helmholtz equation.

Out of six approaches, the DSC-RSK method is the only method that is able to deliver extremely 
accurate solutions at grid density of 2.1 PPW for this problem. It also requires the smallest CPU 
time for the generation of the weights.

\subsection{Helmholtz equation with a multiple frequency source term}\label{helm}

Here we compare the FD and  DSC methods
when there is significant amplitude in a wide range of wavenumbers.
To this end, we consider the Helmholtz equation with a multiple frequency source term
\begin{eqnarray}
u_{xx}(x) + k^2 u(x) & = & f(x), \quad x \in [0,\pi]
\nonumber\\ 
\label{Helmholtz2} 
u(0) & = & u(\pi) = \frac{k}{2}+1, \\
f(x) & = & \sum_{j=0}^{k/2-1} (k^2- 4 j^2) \cos (2 j x), \nonumber \\
u(x) & = & \sum_{j=0}^{k/2}  \cos (2 j x). \nonumber
\end{eqnarray}
By considering $N=526$ and $k=500$, the  grid density of this
problem is also PPW=2.1. The Fourier frequency response of the solution is
depicted in Fig. \ref{fig.Helmsolu} (b). 
This is a problem that involves almost all available spectrum of the discrete
Fourier domain. It is a good example to illustrate ``Proof'''s   claim.

Since this problem requires a large stencil, 
we examine the performance of different methods over a range of large $M$ values, 
from $M=50$ to 500. 
For the Sech method, the near optimized $D=0.28$ is used for all $M$.
For the DSC-RSK method, we choose  $(M,r)=(50,9.0)$, (100,27.9), (150,35.1), (200,38.1),
(250,42.3), (300,46.9), (350,49.6), (400,53.5), (450,56.0), and (500,60.0)

$L_{\infty}$ errors are listed  in Table \ref{table.Helm2}. 
Execution time for generating the first and second derivative weights 
is also reported for each approach. These results are similar to those in 
the previous example. The Euler and MEuler methods show a great deal of improvement 
from their last performance. However, their results are still up to 12 orders
less accurate than those of the DSC-RSK method at certain $M$ values. 

For this multiple frequency problem, it is noted that the accuracy of the FD method does 
not improve much as the $M$ is increased from 50 to 400. It fails to converge  
when $M>400$, indicating the possible instability in the FD method for this class of 
multiple frequency problems. The DSC-RSK method outperform the FD method by up to 
10 orders of magnitude. Out of six methods studied, the DSC-RSK method is the only one
that is able to deliver high accuracy for this multiple frequency problem.

\begin{table}[!tb]       
\caption{Numerical results and CPU times in second 
for solving the  Helmholtz equation with a multiple frequency source term.
Here NaN and --- stand for ``not a number'' and ``not available'', respectively.
}
\label{table.Helm2}      
\begin{center}      
\begin{tabular}{lllllll}
\hline
\hline
 & \multicolumn{2}{c}{FD} & \multicolumn{2}{c}{Boyd's FD} & 
\multicolumn{2}{c}{Sech} \\
\hline
$M$ & Error & CPU  & Error & CPU  & Error & CPU \\
\hline
50 &   $5.06(+0)$ & $8.92(-4)$ & $5.06(+0)$ & $4.88(-4)$ & $1.78(-1)$ & 3.998\\
100 &  $5.48(+0)$ & $3.12(-3)$ & NaN & $3.45(-4)$        & $1.51(-1)$ & 15.73\\
150 &  $2.37(+0)$ & $7.02(-3)$ & --- & ---               & $4.76(-2)$ & 35.38\\
200 &  $4.06(+0)$ & $1.23(-2)$ & --- & ---               & $2.07(-2)$ & 62.98\\
250 &  $1.14(+0)$ & $1.91(-2)$ & --- & ---               & $2.75(-2)$ & 98.56\\
300 &  $6.60(-1)$ & $2.75(-2)$ & --- & ---               & ---  & --- \\
350 &  $4.27(-1)$ & $3.83(-2)$ & --- & ---               & ---  & --- \\
400 &  $2.87(-1)$ & $4.89(-2)$ & --- & ---               & ---  & --- \\
450 &  NaN        & $7.20(-2)$ & --- & ---               & ---  & --- \\
500 &  NaN        & $1.78(-1)$ & --- & ---               & ---  & --- \\
\hline
 & \multicolumn{2}{c}{Euler} & \multicolumn{2}{c}{MEuler} & 
\multicolumn{2}{c}{RSK}\\
\hline
$M$ & Error & CPU  & Error & CPU  & Error & CPU \\
\hline
50 &  $3.56(-1)$ & $4.78(-4)$ & $3.56(-1)$ & $5.00(-4)$  &  $4.63(+0)$ & $3.76(-4)$ \\
100 & $4.74(-1)$ & $4.01(-4)$ & $4.74(-1)$ & $3.92(-4)$  &  $4.76(-3)$ & $2.06(-4)$ \\
150 & $1.80(-1)$ & $5.93(-4)$ & $1.80(-1)$ & $5.95(-4)$  &  $7.79(-5)$ & $3.34(-4)$ \\
200 & $4.61(+1)$ & $5.81(-4)$ & $4.61(+1)$ & $5.78(-4)$  &  $2.52(-7)$ & $4.22(-4)$ \\
250 & $4.61(+1)$ & $6.86(-4)$ & $4.61(+1)$ & $6.83(-4)$  &  $1.17(-8)$ & $5.14(-4)$ \\
300 & $4.61(+1)$ & $2.67(-2)$ & $4.61(+1)$ & $2.68(-2)$  &  $5.36(-10)$ & $6.14(-4)$ \\
350 & $4.61(+1)$ & $2.04(-2)$ & $4.61(+1)$ & $2.05(-2)$  &  $1.87(-11)$ & $7.13(-4)$ \\
400 & $4.61(+1)$ & $2.59(-2)$ & $4.61(+1)$ & $2.63(-2)$  &  $7.32(-12)$ & $8.51(-4)$ \\
450 & $4.61(+1)$ & $2.73(-2)$ & $4.61(+1)$ & $2.83(-2)$  &  $7.31(-12)$ & $9.11(-4)$ \\
500 & $4.61(+1)$ & $3.64(-2)$ & $4.61(+1)$ & $3.67(-2)$  &  $7.39(-12)$ & $1.01(-3)$ \\
\hline
\hline
\end{tabular}       
\end{center}      
\end{table}    

~

\subsection{Unsteady hyperbolic equation with a few wavenumbers}\label{few}

\begin{figure*}[!tb]
\begin{center}
(a) \psfig{figure=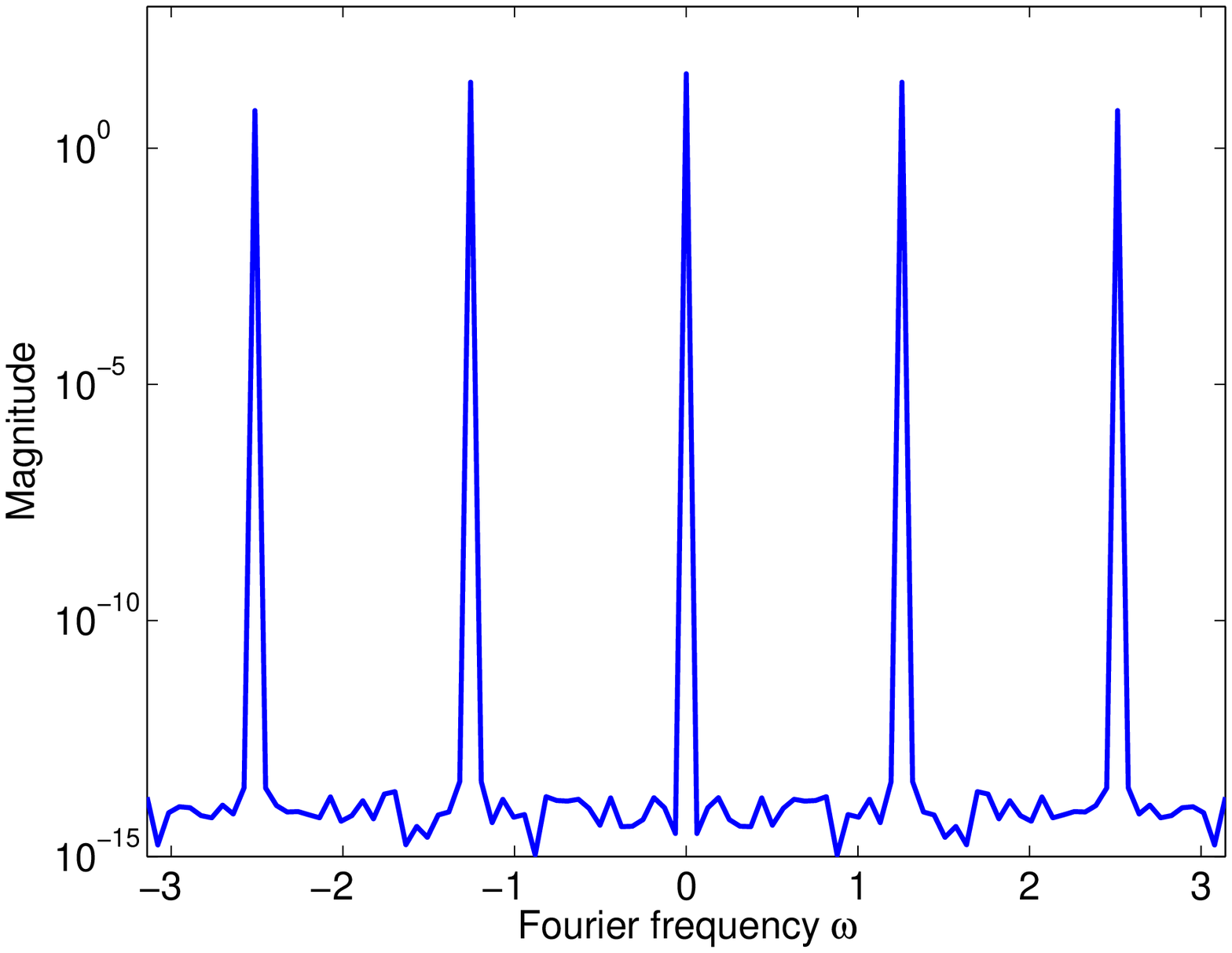,width=0.45\linewidth,height=6.5cm} 
(b) \psfig{figure=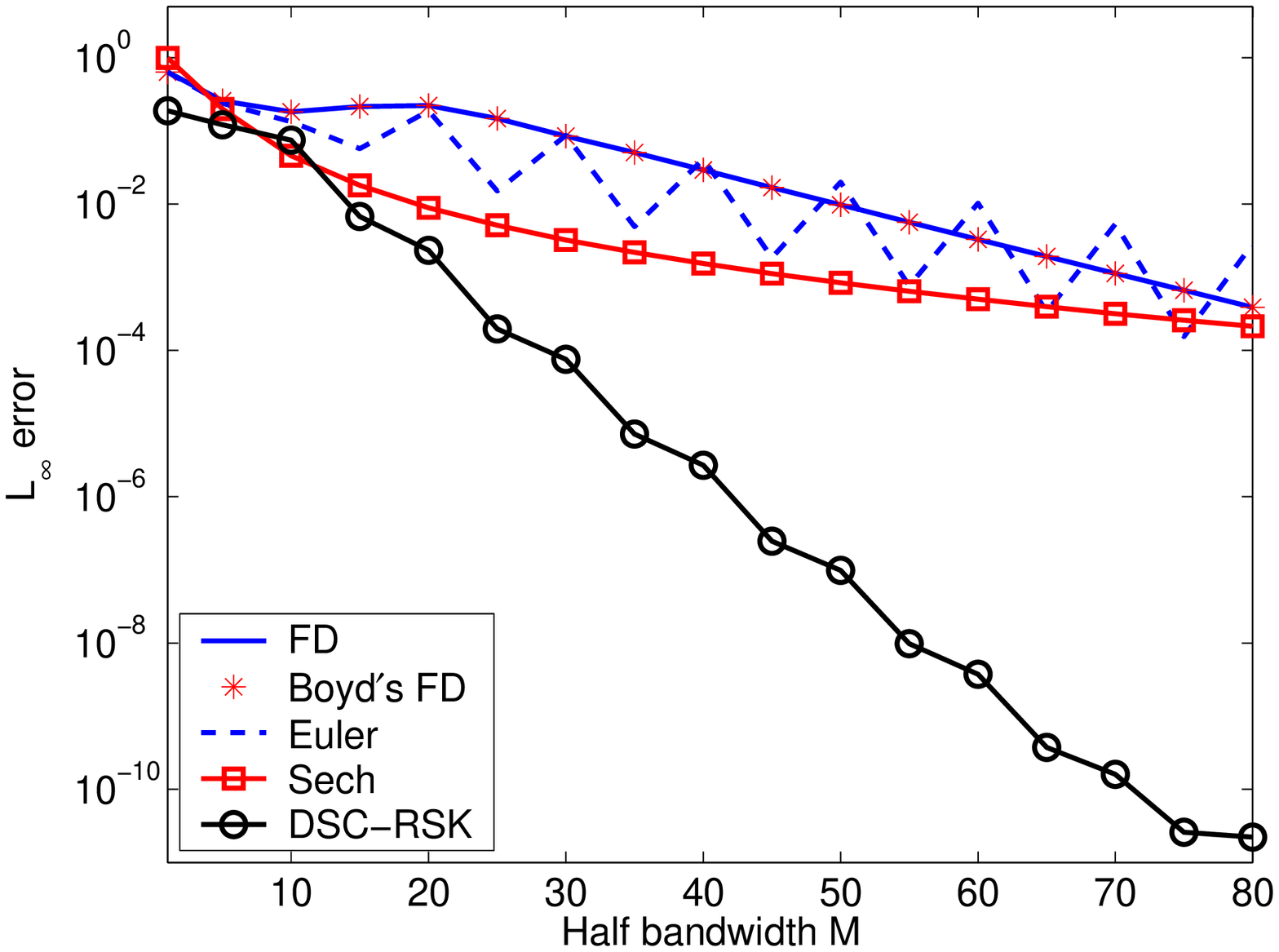,width=0.45\linewidth,height=6.5cm}  \\
\end{center}
\caption{ The unsteady hyperbolic equation with a few wavenumbers.
(a) Fourier frequency response of $u(0,x)$;
(b) Numerical errors.
}
\label{fig.hyperb1} 
\end{figure*}

We consider the time-dependent equation

\begin{eqnarray}\label{hyperbolic}
u_{t} & = & - t^2 u_x,  \quad x \in [-1,1]
\nonumber\\ 
u(0,x) & = & \sin^4 (k \pi x), \\
u(t,x) & = & \sin^4 (k \pi (x-t^3/3)), \nonumber
\end{eqnarray}
with periodical boundary condition. 
Fourier frequency response of $u(0,x)$ is given in Fig. \ref{fig.hyperb1} (a).

We use the classical fourth order Runge-Kutta (RK4) scheme for the time discretization.  
To compare the spatial discretization errors of various higher order approaches, a 
sufficiently small time increment $\Delta t$ has to be used. Here, we choose
$N=101$, $k=10$, $t=1$, and $\Delta t = 5 \times 10^{-5}$. 
It was numerically tested that a smaller $\Delta t$ would not yield significantly 
more accurate results. In the Sech method, a near optimal parameter $D=0.36$ is used.
In the DSC-RSK method, we choose $(M,r)=(1,1.9)$, (5,2.4), (10,4.1), (15,5.3), (20,5.9),
(25,6.7), (30,7.2), (35,7.8), (40,8.2), (45,8.8), (50,9.2), (55,9.7),
(60,10.0), (65,10.5), (70,10.8), (75,11.4), and (80,11.5).

The $L_{\infty}$ errors are depicted in Fig. \ref{fig.hyperb1}(b). 
Boyd's FD method yields the same results as the FD method. 
It is clear that for this highly oscillatory solution, the DSC-RSK method is
the best for all $M$ values examined, including small stencils, like $M=1$, 5 and 10. 
When $M > 75$, the DSC-RSK error does not decay any more, because  the accuracy 
limit of the RK4 temporal integration is reached. The accuracies of other 
approaches are far away from the limit. It is seen that the DSC-RSK method is over
$10^{6}$ times more accurate than all other approaches examined for some large $M$.

\subsection{Unsteady hyperbolic equation with a  wave comb}\label{comb}

\begin{figure*}[!tb]
\begin{center}
a) \psfig{figure=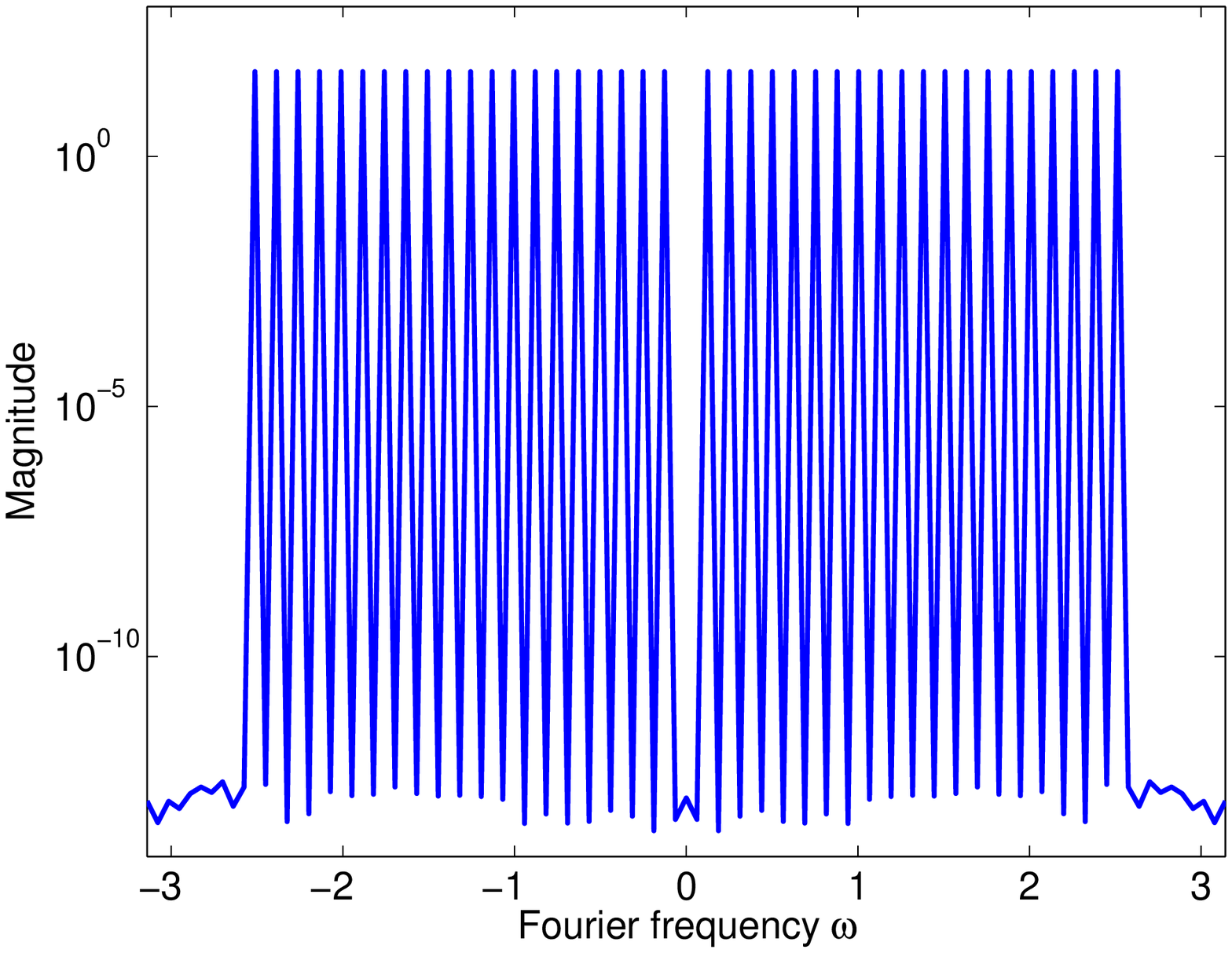,width=0.45\linewidth,height=6.5cm}
b) \psfig{figure=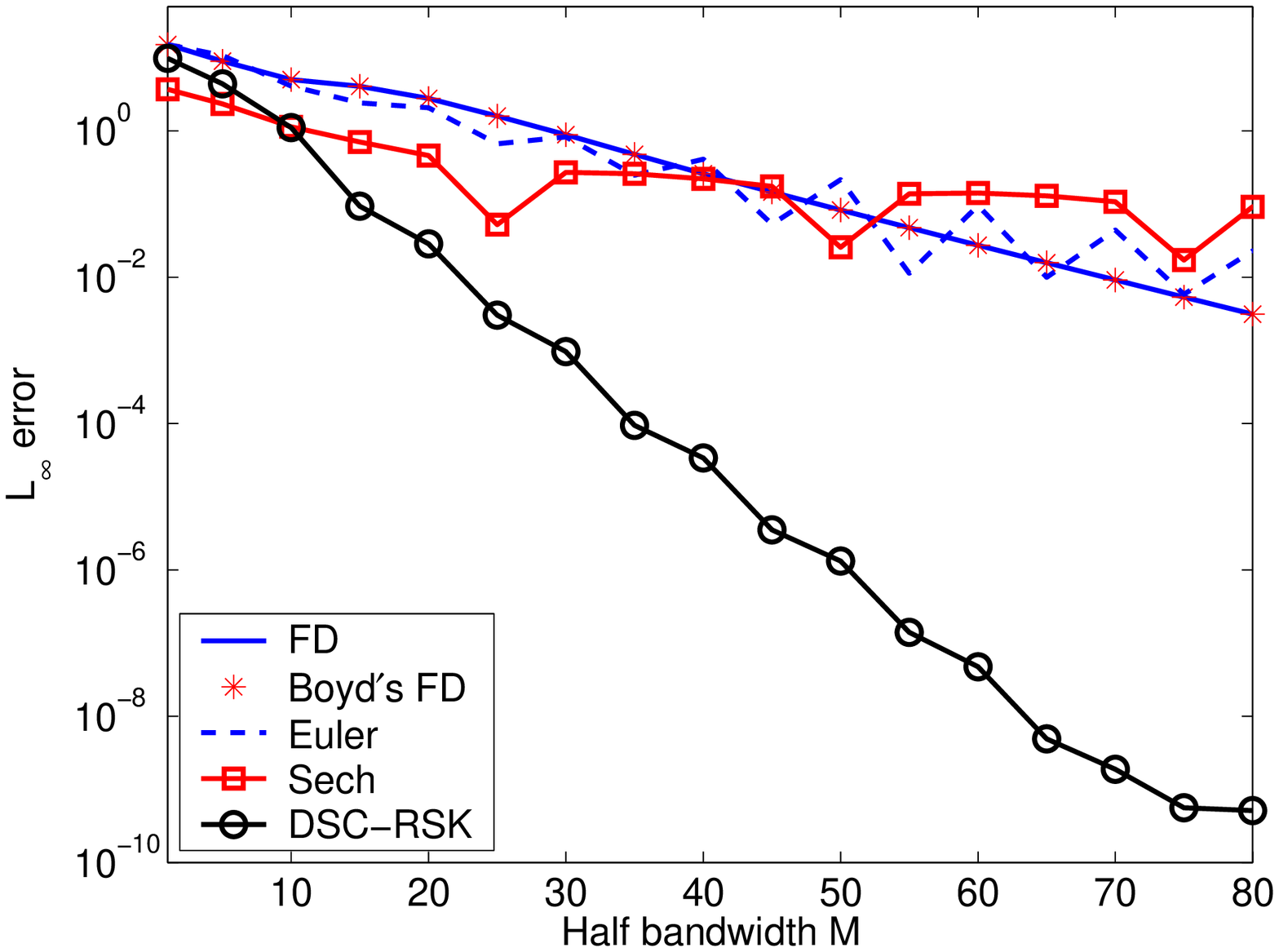,width=0.45\linewidth,height=6.5cm} \\
\end{center}
\caption{ The unsteady hyperbolic equation with a wave comb.
(a) Fourier frequency response;
(b) Numerical errors.
}
\label{fig.hyperb2} 
\end{figure*}

We next consider Eq. (\ref{hyperbolic}) with a wide-band frequency solution
\begin{equation}
u(t,x)= \sum_{j=1}^{k/2} \sin (2 j \pi (x-t^3/3)).
\end{equation}
Fig. \ref{fig.hyperb2} (a) depicts the Fourier frequency response, which 
has a wide range of  wavenumbers. To examine ``Proof'''s   claim about the FD method, 
we have chosen the amplitude in small wavenumbers to be as large as those in large 
wavenumbers.

Here we choose $N=101$, $k=40$, and  $\Delta t = 5 \times 10^{-5}$ in our computation.
The free parameters $D$ and $r$ remain the same as in the last case. 
The results at $t=1$ are shown in Fig. \ref{fig.hyperb2} (b). Although none of the methods 
performs as well as in the previous case, the Sech method has the most dramatical accuracy 
reduction. Its error has increased 1000 times, indicating  its sensitivity to the frequency 
distribution.

The DSC-RSK method keeps its edge over all other methods examined --- it is up to a million 
times more accurate than the FD and Sech methods at large stencils. In fact, it outperforms 
the FD method over all stencil widths examined.

\subsection{Two dimensional Navier-Stokes equation}\label{nav}

We next consider the two dimensional (2D) Navier-Stokes equation
in its primitive variable, describing an incompressible fluid flow
\begin{eqnarray}
\frac{\partial {\bf u}}{\partial t} + {\bf u} \cdot \nabla {\bf u}
& = & -\nabla p + \frac{1}{\mbox{Re}}  \nabla^2 {\bf u}, \\
\nabla \cdot {\bf u} & = & 0,
\end{eqnarray}
where ${\bf u}=(u,v)$ is the velocity vector, $p$ is the pressure
and Re is the Reynolds number.
In a 2D square domain $[0, 2\pi] \times [0, 2\pi]$ with periodic
boundary conditions, the initial values are taken as
\begin{eqnarray}
u(0,x,y) & = & - \cos(kx) \sin(ky), \\
v(0,x,y) & = & \sin(kx) \cos(ky). \nonumber
\end{eqnarray}
The analytical solution is
\begin{eqnarray}
u(t,x,y) & = & - \cos(kx) \sin(ky) \exp ( -2k^2 t/ \mbox{Re}), \nonumber \\
v(t,x,y) & = & \sin(kx) \cos(ky)\exp ( -2k^2 t/ \mbox{Re}), \\
p(t,x,y) & = & -\frac{1}{4} [ \cos(2kx) + \cos(2ky)]\exp ( -4k^2 t/ \mbox{Re}).\nonumber 
\end{eqnarray}

The present study has two interesting points in comparing with the previous studies,
i.e., it is multi-dimensional and it involves nonlinearity.
We adopt the Adams-Bashforth-Crank-Nicolson (ABCN) scheme for the time 
discretization and treatment of pressure
\begin{eqnarray}
\frac{1}{2\mbox{Re}}  \nabla^2 {\bf u}^{n+1} - \frac{1}{\Delta t} {\bf u}^{n+1}
& = & \nabla p^{n+1/2} + {\bf S}^{n} \\
\nabla^2 p^{n+1/2} + \nabla \cdot {\bf S}^{n}  & = & 0, \nonumber
\end{eqnarray}
where the source vector ${\bf S}^{n}$ is given by
\begin{equation}
{\bf S}^{n}=-\frac{{\bf u}^{n}}{\Delta t} + 
\frac{3({\bf u}^{n} \cdot \nabla) {\bf u}^{n}-({\bf u}^{n-1} \cdot \nabla){\bf u}^{n-1}}{2}
- \frac{1}{2\mbox{Re}}  \nabla^2 {\bf u}^{n}.
\end{equation}

\begin{figure*}[!tb]
\begin{center}
\psfig{figure=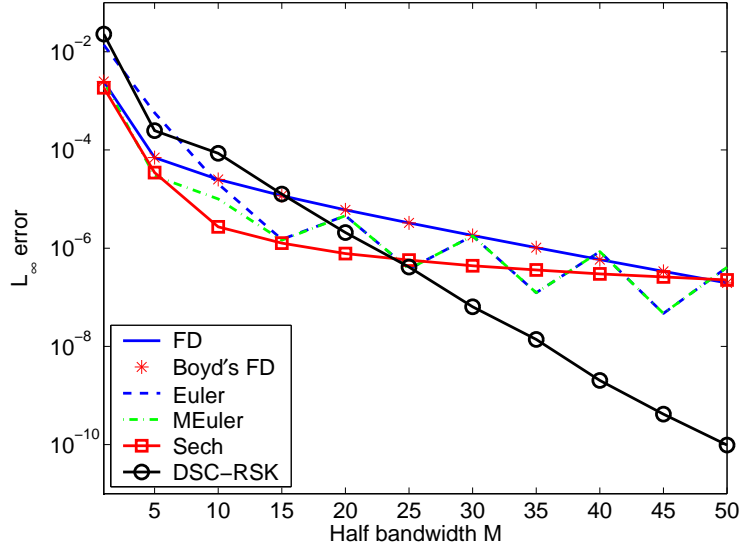,width=0.6\linewidth} 
\end{center}
\caption{Numerical errors of the Navier-Stokes equation.}
\label{fig.ns} 
\end{figure*}

We consider a high frequency problem with $k=10$ and $\mbox{Re}=100$.
By using 51 nodes along each direction, the grid density for solving the Poisson equation of the 
pressure is PPW=2.5, while that for the velocity field is PPW=5. 
In order to compare spatial discretization errors, 
a sufficiently small time increment, $\Delta t=1.0 \times 10^{-5}$, is used.
For the Sech method, a near optimal $D=0.18$ is used. For the DSC-RSK method, we 
choose $(M,r)=(1,0.8)$, (5,1.7), (10,4.6), (15,5.4), (20,6.1),
(25,6.8), (30,7.3), (35,7.9), (40,8.3), (45,8.8), and (50,9.3).
The $L_{\infty}$ errors in $u$ after 1000 time steps are depicted in Fig. \ref{fig.ns}. 
It is clear that the DSC-RSK method is more accurate than other approaches for this high
frequency problem when $M > 25$. At $M=50$, the DSC-RSK method is over
a thousand times more accurate than the Sech and FD methods.

\subsection{Quantum eigenvalue problem}\label{quantum}

Finally, we consider a quantum eigenvalue problem
\begin{equation}\label{harmonic}
\left[ -\frac{1}{2} \frac{d^2}{d x^2} + V(x) \right]
\Phi_n = E_n \Phi_n,
\end{equation}
with harmonic oscillator potential $V(x)=x^2/2$. The eigenfunctions
$\Phi_n$ are the standard Hermite function series. The eigenvalues 
are $E_n=n+\frac{1}{2}$, where $n=0,1,2,\cdots$. 

\setlongtables
\begin{longtable}{llllllll}
\caption{Numerical errors in the quantum eigenvalue analysis with $M=50$ and $N=51$.} 
\label{table.ho50}  \\
\hline
\hline
Mode & Sinc & DSC-RSK & FD & Boyd's FD & Euler & MEuler & Sech \\
\hline
0  & $1.40(-14)$ & $9.10(-15)$ & $1.38(-13)$ & $7.17(-13)$ & $2.67(-15)$ & $2.67(-15)$ & $3.46(-14)$ \\ 
1  & $1.48(-14)$ & $3.85(-15)$ & $3.75(-14)$ & $2.23(-13)$ & $4.00(-15)$ & $4.00(-15)$ & $2.75(-14)$ \\ 
2  & $5.33(-15)$ & $5.33(-16)$ & $2.56(-14)$ & $1.18(-13)$ & $5.51(-15)$ & $5.51(-15)$ & $3.30(-14)$ \\ 
3  & $3.04(-15)$ & $1.27(-16)$ & $1.28(-14)$ & $7.74(-14)$ & $7.99(-15)$ & $7.99(-15)$ & $4.15(-14)$ \\ 
4  & $5.92(-16)$ & $5.92(-16)$ & $1.12(-14)$ & $5.51(-14)$ & $7.11(-15)$ & $7.11(-15)$ & $5.01(-14)$ \\ 
5  & $9.69(-16)$ & $4.84(-16)$ & $8.07(-15)$ & $4.00(-14)$ & $6.94(-15)$ & $6.94(-15)$ & $1.26(-14)$ \\ 
6  & $2.73(-16)$ & $1.64(-15)$ & $4.24(-15)$ & $3.07(-14)$ & $7.52(-15)$ & $7.52(-15)$ & $9.34(-13)$ \\ 
7  & $1.18(-15)$ & $4.74(-16)$ & $6.99(-15)$ & $2.21(-14)$ & $5.45(-15)$ & $5.45(-15)$ & $1.04(-11)$ \\ 
8  & $8.36(-16)$ & $2.09(-16)$ & $1.94(-14)$ & $4.60(-15)$ & $2.09(-15)$ & $2.09(-15)$ & $9.17(-11)$ \\ 
9  & $7.48(-16)$ & $0.00$ & $1.01(-13)$ & $8.23(-14)$ & $2.43(-15)$ & $2.43(-15)$ & $7.03(-10)$ \\ 
10  & $5.07(-16)$ & $0.00$ & $5.60(-13)$ & $5.45(-13)$ & $1.51(-13)$ & $1.51(-13)$ & $4.76(-9)$ \\ 
11  & $4.63(-16)$ & $1.54(-16)$ & $2.87(-12)$ & $2.86(-12)$ & $7.81(-13)$ & $7.81(-13)$ & $2.86(-8)$ \\ 
12  & $7.11(-16)$ & $4.26(-16)$ & $1.35(-11)$ & $1.35(-11)$ & $1.58(-12)$ & $1.58(-12)$ & $1.53(-7)$ \\ 
13  & $1.32(-16)$ & $2.63(-16)$ & $5.83(-11)$ & $5.83(-11)$ & $2.09(-12)$ & $2.09(-12)$ & $7.37(-7)$ \\ 
14  & $3.68(-16)$ & $1.22(-16)$ & $2.33(-10)$ & $2.33(-10)$ & $2.31(-11)$ & $2.31(-11)$ & $3.18(-6)$ \\ 
15  & $1.38(-15)$ & $1.15(-15)$ & $8.68(-10)$ & $8.68(-10)$ & $6.40(-11)$ & $6.40(-11)$ & $1.23(-5)$ \\ 
16  & $6.46(-16)$ & $1.29(-15)$ & $3.01(-9)$ & $3.01(-9)$ & $3.65(-11)$ & $3.65(-11)$ & $4.31(-5)$ \\ 
17  & $7.31(-14)$ & $7.31(-14)$ & $9.79(-9)$ & $9.79(-9)$ & $3.23(-10)$ & $3.23(-10)$ & $1.35(-4)$ \\ 
18  & $5.57(-14)$ & $6.32(-14)$ & $2.99(-8)$ & $2.99(-8)$ & $1.19(-9)$ & $1.19(-9)$ & $3.78(-4)$ \\ 
19  & $2.87(-12)$ & $2.88(-12)$ & $8.61(-8)$ & $8.61(-8)$ & $1.54(-9)$ & $1.54(-9)$ & $9.48(-4)$ \\ 
20  & $1.78(-12)$ & $1.98(-12)$ & $2.34(-7)$ & $2.34(-7)$ & $1.96(-9)$ & $1.96(-9)$ & $2.12(-3)$ \\ 
21  & $8.66(-11)$ & $8.70(-11)$ & $6.04(-7)$ & $6.04(-7)$ & $1.24(-8)$ & $1.24(-8)$ & $4.26(-3)$ \\ 
22  & $4.42(-11)$ & $5.04(-11)$ & $1.48(-6)$ & $1.48(-6)$ & $2.17(-8)$ & $2.17(-8)$ & $7.71(-3)$ \\ 
23  & $1.99(-9)$ & $1.99(-9)$ & $3.45(-6)$ & $3.45(-6)$ & $3.79(-9)$ & $3.79(-9)$ & $1.27(-2)$ \\ 
24  & $7.87(-10)$ & $9.00(-10)$ & $7.67(-6)$ & $7.67(-6)$ & $7.55(-8)$ & $7.55(-8)$ & $1.94(-2)$ \\ 
25  & $3.57(-8)$ & $3.59(-8)$ & $1.63(-5)$ & $1.63(-5)$ & $1.40(-7)$ & $1.40(-7)$ & $2.78(-2)$ \\ 
26  & $1.12(-8)$ & $1.33(-8)$ & $3.33(-5)$ & $3.33(-5)$ & $1.23(-7)$ & $1.23(-7)$ & $3.77(-2)$ \\ 
27  & $4.93(-7)$ & $4.94(-7)$ & $6.44(-5)$ & $6.44(-5)$ & $7.59(-7)$ & $7.59(-7)$ & $4.90(-2)$ \\ 
28  & $1.05(-7)$ & $1.29(-7)$ & $1.23(-4)$ & $1.23(-4)$ & $1.08(-6)$ & $1.08(-6)$ & $6.16(-2)$ \\ 
29  & $5.34(-6)$ & $5.36(-6)$ & $2.15(-4)$ & $2.15(-4)$ & $6.43(-6)$ & $6.43(-6)$ & $7.52(-2)$ \\ 
30  & $7.32(-7)$ & $1.01(-6)$ & $4.04(-4)$ & $4.04(-4)$ & $2.05(-8)$ & $2.05(-8)$ & $8.98(-2)$ \\ 
31  & $4.40(-5)$ & $4.41(-5)$ & $6.01(-4)$ & $6.01(-4)$ & $4.02(-5)$ & $4.02(-5)$ & $1.05(-1)$ \\ 
32  & $7.74(-7)$ & $2.83(-6)$ & $1.24(-3)$ & $1.24(-3)$ & $2.88(-6)$ & $2.88(-6)$ & $1.21(-1)$ \\ 
33  & $2.81(-4)$ & $2.82(-4)$ & $1.33(-3)$ & $1.33(-3)$ & $2.81(-4)$ & $2.81(-4)$ & $1.38(-1)$ \\ 
34  & $2.87(-5)$ & $1.35(-5)$ & $3.86(-3)$ & $3.86(-3)$ & $2.54(-5)$ & $2.54(-5)$ & $1.56(-1)$ \\ 
35  & $1.33(-3)$ & $1.33(-3)$ & $1.93(-3)$ & $1.93(-3)$ & $1.34(-3)$ & $1.34(-3)$ & $1.74(-1)$ \\ 
36  & $3.82(-4)$ & $3.07(-4)$ & $1.16(-2)$ & $1.16(-2)$ & $3.62(-4)$ & $3.62(-4)$ & $1.92(-1)$ \\ 
37  & $4.83(-3)$ & $4.84(-3)$ & $4.36(-4)$ & $4.36(-4)$ & $4.83(-3)$ & $4.83(-3)$ & $2.11(-1)$ \\ 
38  & $1.48(-3)$ & $1.19(-3)$ & $2.04(-2)$ & $2.04(-2)$ & $1.54(-3)$ & $1.54(-3)$ & $2.31(-1)$ \\ 
39  & $1.34(-2)$ & $1.34(-2)$ & $6.66(-3)$ & $6.66(-3)$ & $1.33(-2)$ & $1.33(-2)$ & $2.50(-1)$ \\ 
40  & $2.95(-4)$ & $1.68(-4)$ & $1.71(-2)$ & $1.71(-2)$ & $1.86(-4)$ & $1.86(-4)$ & $2.70(-1)$ \\ 
41  & $3.06(-2)$ & $3.06(-2)$ & $2.29(-2)$ & $2.29(-2)$ & $3.07(-2)$ & $3.07(-2)$ & $2.91(-1)$ \\ 
42  & $1.19(-2)$ & $1.24(-2)$ & $1.05(-3)$ & $1.05(-3)$ & $1.19(-2)$ & $1.19(-2)$ & $3.11(-1)$ \\ 
43  & $5.93(-2)$ & $5.93(-2)$ & $5.13(-2)$ & $5.13(-2)$ & $5.95(-2)$ & $5.95(-2)$ & $3.32(-1)$ \\ 
44  & $3.82(-2)$ & $3.86(-2)$ & $2.77(-2)$ & $2.77(-2)$ & $3.81(-2)$ & $3.81(-2)$ & $3.53(-1)$ \\ 
45  & $1.05(-1)$ & $1.06(-1)$ & $9.61(-2)$ & $9.61(-2)$ & $1.06(-1)$ & $1.06(-1)$ & $3.76(-1)$ \\ 
46  & $8.35(-2)$ & $8.39(-2)$ & $7.26(-2)$ & $7.26(-2)$ & $8.30(-2)$ & $8.30(-2)$ & $3.91(-1)$ \\ 
47  & $1.76(-1)$ & $1.76(-1)$ & $1.67(-1)$ & $1.67(-1)$ & $1.76(-1)$ & $1.76(-1)$ & $4.31(-1)$ \\ 
48  & $1.53(-1)$ & $1.53(-1)$ & $1.43(-1)$ & $1.43(-1)$ & $1.52(-1)$ & $1.52(-1)$ & $4.21(-1)$ \\ 
49  & $3.14(-1)$ & $3.14(-1)$ & $2.93(-1)$ & $2.93(-1)$ & $3.15(-1)$ & $3.15(-1)$ & $5.49(-1)$ \\ 
50  & $2.89(-1)$ & $2.90(-1)$ & $2.68(-1)$ & $2.68(-1)$ & $2.88(-1)$ & $2.88(-1)$ & $5.25(-1)$ \\ 
\hline
\hline
\end{longtable}

\setlongtables
\begin{longtable}{lllllll}
\caption{Numerical errors in the quantum eigenvalue analysis with $M=200$ and $N=201$.} 
\label{table.ho200}  \\
\hline
\hline
Mode & Sinc & DSC-RSK & FD & Euler & MEuler & Sech \\
\hline
0  & $8.49(-13)$ & $7.03(-13)$ & $8.70(-13)$ & $1.00$ & $1.00$ & $1.33(-13)$ \\
4  & $2.11(-14)$ & $6.12(-15)$ & $7.80(-14)$ & $9.86(-1)$ & $9.86(-1)$ & $1.42(-14)$ \\
8  & $2.28(-14)$ & $6.48(-15)$ & $9.57(-14)$ & $9.71(-1)$ & $9.71(-1)$ & $1.02(-14)$ \\
12  & $1.21(-14)$ & $1.02(-14)$ & $5.74(-14)$ & $9.55(-1)$ & $9.55(-1)$ & $9.66(-14)$ \\
16  & $3.66(-15)$ & $4.52(-15)$ & $5.23(-14)$ & $9.40(-1)$ & $9.40(-1)$ & $1.07(-13)$ \\
20  & $3.81(-15)$ & $9.70(-15)$ & $4.75(-14)$ & $9.24(-1)$ & $9.24(-1)$ & $1.52(-13)$ \\
24  & $3.48(-15)$ & $3.04(-15)$ & $4.21(-14)$ & $9.09(-1)$ & $9.09(-1)$ & $1.84(-13)$ \\
28  & $3.74(-16)$ & $4.24(-15)$ & $3.54(-14)$ & $8.93(-1)$ & $8.93(-1)$ & $1.36(-13)$ \\
32  & $4.15(-15)$ & $3.94(-15)$ & $3.24(-14)$ & $8.78(-1)$ & $8.78(-1)$ & $1.24(-13)$ \\
36  & $2.53(-15)$ & $2.92(-15)$ & $2.94(-14)$ & $8.63(-1)$ & $8.63(-1)$ & $8.10(-14)$ \\
40  & $7.02(-16)$ & $1.75(-15)$ & $2.72(-14)$ & $8.47(-1)$ & $8.47(-1)$ & $2.25(-13)$ \\
44  & $4.79(-16)$ & $4.79(-16)$ & $2.44(-14)$ & $8.32(-1)$ & $8.32(-1)$ & $2.32(-13)$ \\
48  & $1.46(-16)$ & $1.03(-15)$ & $2.27(-14)$ & $8.16(-1)$ & $8.16(-1)$ & $5.51(-13)$ \\
52  & $0.00$ & $1.35(-16)$ & $2.18(-14)$ & $8.01(-1)$ & $8.01(-1)$ & $1.18(-12)$ \\
56  & $1.26(-16)$ & $2.51(-16)$ & $1.96(-14)$ & $7.85(-1)$ & $7.85(-1)$ & $1.53(-12)$ \\
60  & $3.52(-16)$ & $2.35(-16)$ & $1.77(-14)$ & $7.70(-1)$ & $7.70(-1)$ & $1.72(-12)$ \\
64  & $0.00$ & $8.81(-16)$ & $1.74(-14)$ & $7.54(-1)$ & $7.54(-1)$ & $1.06(-12)$ \\
68  & $2.07(-16)$ & $0.00$ & $1.53(-14)$ & $7.39(-1)$ & $7.39(-1)$ & $9.25(-14)$ \\
72  & $7.84(-16)$ & $3.92(-16)$ & $1.41(-14)$ & $7.23(-1)$ & $7.23(-1)$ & $1.01(-12)$ \\
76  & $5.57(-16)$ & $3.72(-16)$ & $1.34(-14)$ & $7.08(-1)$ & $7.08(-1)$ & $1.45(-12)$ \\
80  & $0.00$ & $0.00$ & $1.29(-14)$ & $6.92(-1)$ & $6.92(-1)$ & $3.57(-13)$ \\
84  & $0.00$ & $3.36(-16)$ & $1.21(-14)$ & $6.77(-1)$ & $6.77(-1)$ & $7.92(-13)$ \\
88  & $0.00$ & $3.21(-16)$ & $1.19(-14)$ & $6.61(-1)$ & $6.61(-1)$ & $1.43(-12)$ \\
92  & $3.07(-16)$ & $3.07(-16)$ & $1.69(-14)$ & $6.46(-1)$ & $6.46(-1)$ & $3.58(-12)$ \\
96  & $2.95(-16)$ & $0.00$ & $1.41(-13)$ & $6.30(-1)$ & $6.30(-1)$ & $4.46(-12)$ \\
100  & $5.66(-16)$ & $2.83(-16)$ & $2.12(-12)$ & $6.15(-1)$ & $6.15(-1)$ & $8.07(-13)$ \\
104  & $5.44(-16)$ & $0.00$ & $2.82(-11)$ & $5.99(-1)$ & $5.99(-1)$ & $5.98(-12)$ \\
108  & $6.55(-16)$ & $6.55(-16)$ & $3.11(-10)$ & $5.84(-1)$ & $5.84(-1)$ & $2.65(-11)$ \\
112  & $6.32(-16)$ & $1.01(-15)$ & $2.86(-9)$ & $5.68(-1)$ & $5.68(-1)$ & $5.61(-11)$ \\
116  & $6.10(-16)$ & $1.22(-16)$ & $2.20(-8)$ & $5.53(-1)$ & $5.53(-1)$ & $8.50(-11)$ \\
120  & $2.19(-14)$ & $2.17(-14)$ & $1.43(-7)$ & $5.37(-1)$ & $5.37(-1)$ & $1.47(-10)$ \\
124  & $1.25(-12)$ & $1.26(-12)$ & $7.88(-7)$ & $5.22(-1)$ & $5.22(-1)$ & $2.11(-10)$ \\
128  & $5.35(-11)$ & $5.39(-11)$ & $3.71(-6)$ & $5.06(-1)$ & $5.06(-1)$ & $1.69(-10)$ \\
132  & $1.70(-9)$ & $1.72(-9)$ & $1.51(-5)$ & $4.91(-1)$ & $4.91(-1)$ & $1.45(-9)$ \\
136  & $3.94(-8)$ & $3.99(-8)$ & $5.36(-5)$ & $4.75(-1)$ & $4.75(-1)$ & $3.91(-8)$ \\
140  & $6.47(-7)$ & $6.57(-7)$ & $1.70(-4)$ & $4.60(-1)$ & $4.60(-1)$ & $6.47(-7)$ \\
144  & $7.14(-6)$ & $7.30(-6)$ & $4.99(-4)$ & $4.44(-1)$ & $4.44(-1)$ & $7.14(-6)$ \\
148  & $4.66(-5)$ & $4.86(-5)$ & $1.60(-3)$ & $4.29(-1)$ & $4.29(-1)$ & $4.66(-5)$ \\
152  & $1.21(-4)$ & $1.39(-4)$ & $5.24(-3)$ & $4.13(-1)$ & $4.13(-1)$ & $1.21(-4)$ \\
156  & $7.86(-5)$ & $1.60(-4)$ & $6.37(-3)$ & $3.98(-1)$ & $3.98(-1)$ & $7.86(-5)$ \\
160  & $1.97(-3)$ & $2.08(-3)$ & $2.57(-3)$ & $3.82(-1)$ & $3.82(-1)$ & $1.97(-3)$ \\
164  & $8.39(-3)$ & $8.48(-3)$ & $4.95(-3)$ & $3.67(-1)$ & $3.67(-1)$ & $8.39(-3)$ \\
168  & $1.91(-2)$ & $1.92(-2)$ & $1.61(-2)$ & $3.51(-1)$ & $3.51(-1)$ & $1.91(-2)$ \\
172  & $3.39(-2)$ & $3.40(-2)$ & $3.12(-2)$ & $3.36(-1)$ & $3.36(-1)$ & $3.39(-2)$ \\
176  & $5.31(-2)$ & $5.32(-2)$ & $5.05(-2)$ & $3.21(-1)$ & $3.21(-1)$ & $5.31(-2)$ \\
180  & $7.74(-2)$ & $7.75(-2)$ & $7.48(-2)$ & $3.05(-1)$ & $3.05(-1)$ & $7.74(-2)$ \\
184  & $1.08(-1)$ & $1.08(-1)$ & $1.05(-1)$ & $2.89(-1)$ & $2.89(-1)$ & $1.08(-1)$ \\
188  & $1.47(-1)$ & $1.47(-1)$ & $1.44(-1)$ & $2.74(-1)$ & $2.74(-1)$ & $1.47(-1)$ \\
192  & $1.99(-1)$ & $1.99(-1)$ & $1.96(-1)$ & $2.59(-1)$ & $2.59(-1)$ & $1.99(-1)$ \\
196  & $2.74(-1)$ & $2.74(-1)$ & $2.69(-1)$ & $2.43(-1)$ & $2.43(-1)$ & $2.74(-1)$ \\
200  & $4.15(-1)$ & $4.15(-1)$ & $4.05(-1)$ & $2.28(-1)$ & $2.28(-1)$ & $4.15(-1)$ \\
\hline
\hline
\end{longtable}

For the Sech method, a near optimal $D=0.28$ is used. In the DSC-RSK method, relatively large values of $r$ are
chosen: $(M,r)= (50,35)$, and (200,90). For $N=51$, and 201, the computational 
domain is chosen as $x \in [-8.7, 8.7]$,  and [-17.6, 17.6], respectively.

This example is chosen to show the subtleness in numerical analysis. A method that is 
good for one problem might not be good for other problems. Detailed analysis is required
for a given problem. In particular, we show that the analysis and claim in 
``Proof''  about the Sinc pseudospectral method is invalid. 

The relative errors in estimating eigenmodes for different $M$
with $N=M+1$ are given in Table \ref{table.ho50} and Table \ref{table.ho200}. 
Since the $O(1/x)$ decay in the Sinc kernel is outpaced by the Gaussian decay of  
the Hermite functions when $|x|$ is large, there is no need for sum acceleration.  
In fact, the Euler-accelerated sinc method (Euler) does not work as well as the 
Sinc does. Clearly, the Sinc method is the best for all $M$ and all eigenmodes.
As it is known that when $r \to \infty$, the DSC-RSK method approaches to the Sinc method,
a very large $r$ is chosen for each $M$, such that the DSC-RSK method results are
very similar to those of the Sinc. Thus, overall, the DSC-RSK  method performs better than the 
FD, Euler, MEuler, and Sech methods.  

What comes as a surprise is the Sech method, which  is able to deliver results of 
double precision accuracy --- in contrast to the fact that it does not reach the accuracy 
of $10^{-7}$ in all other test cases. 
Overall, although it is  not as accurate as the Sinc and DSC-RSK method, the Sech method performs 
better than the FD, Euler and MEuler methods.

When $M \ge 50$, the results of the Euler and MEuler  methods are essentially the same. 
In particular, for $M=200$, they all yield wrong results.

\section{Analysis of ``Proof'''s miscellaneous claims}\label{further}

In this section, we further analyze ``Proof'''s claims. For convenience, these
claims will be presented in an itemized style. 

\begin{itemize}
\item {\sc [``Proof''- Section 1]} 
``It is impossible to review the work of this very
prolific author in detail, but of the more than one hundred articles listed at 
http://www.math.msu.edu/wei/, most use the DSC or LDAF schemes.''

{\sc [Analysis]}
The contrary is true: More than half of the articles do not use  the DSC or LDAF schemes.

\item {\sc [``Proof''- Section 1]}
``Since the finite difference weights are {\it nearly} Gaussian, one cannot escape the conclusion that the
 LDAF/DSC methods are really just high order finite difference methods in disguise!''

{\sc [Analysis]}
There is no disguise about the relation between the DSC and high order finite differences
in our work, as the relation  has been spelled out many times in the DSC literature 
(e.g. in  \cite{weijpa,weiijnme}). It is surprising to see this accusation as the author 
was informed of the existence of these references in a comment on his draft.
Note that all collocation schemes, including the Daubechies wavelet collocation method,
can  be cast into the finite difference form \cite{Fornberg1}.
However, the DSC method has also been formulated in the Galerkin form \cite{Lim}
for solving differential equations, which is not a high order finite difference method.

\item {\sc [``Proof''- Section 1]}
``The LDAF/DSC is also a special case of Boyd's earlier theory of sum-accelerated 
pseudospectral methods: special in that the weighting function is a Gaussian.'' 

{\sc [Analysis]}
This is an unfounded claim.
First, as a computational (discrete) realization of singular convolution kernels,
including kernels of delta type, Abel type and Hilbert type \cite{weijcp99}, 
the DSC has little to do with Boyd's sum-accelerated method. 

Second, the DSC-RSK approximation of derivatives does not fit into Boyd's sum-acceleration 
form $\frac{{\rm d} u}{{\rm d} x} \approx \sum_{j=-M}^{M} w^{\rm DSC}_{Mj}d_j^{(1),{\rm sinc}}$ 
in general --- there is an extra term, see Section \ref{methods}. Only in a very special
case, the on-grid differentiation for the first derivative, can the DSC-RSK  method reduce to a form
similar to Boyd's sum-accelerated pseudospectral methods.  However, off-grid differentiations,
which have been widely used in staggered grids in the previous DSC-RSK applications, cannot.
More extra terms occur even for on-grid differentiations when higher order derivatives are 
considered.  

Finally, the DSC-RSK method is constructed by accelerating the Whittaker-Shannon-Kotel'nikov sampling 
(interpolation) kernel. 

\item {\sc [``Proof''- Section 2]}
``Fourier analysis has been widely used to analyze difference formulas ever since this was popularized by 
von Neuman(n). The reason is that the Fourier basis function, $\exp^{(ikx)}$, is an eigenfunction of both 
the differentiation operator and also of all possible difference formulas. This implies that the 
accuracy of difference formulas can be assessed - and improved - merely by comparing the eigenvalues.''

{\sc [Analysis]}

First, it is  inappropriate to state that ``the accuracy of difference formulas can be assessed - and improved 
- merely by comparing the eigenvalues''. In general, the performance of a difference formula for approximating 
a derivative of a function cannot be assessed by its performance on differentiating $\exp^{(ikx)}$. The reason 
is the follows. In differentiating $\exp^{(ikx)}$, the (discrete) Fourier method is exact and all difference 
formulas are less accurate. This does not justify the Fourier method to be a universal standard for
judging a numerical method's performance of differentiation of a general function, because the 
performance always depends on the functional class, regularity, decay rate, etc. A simple counterexample 
to ``Proof'''s   logic is that in differentiating an $n$th order polynomial, an $n$th order FD method is 
{\it exact}, while the Fourier method is {\it not}.   

Second, in solving differential equations, the superiority of differentiating $\exp^{(ikx)}$ does not 
establish the superiority of a numerical method in general, since differential equations 
might involve functions of different mathematical properties, and be either nonlinear, or 
with variable coefficients, or with other complex non-differentiating terms. For example, the  
Hermite spectral method does not give an exact result for the differentiation of $\exp^{(ikx)}$ 
as the Fourier method does. However, it can give the exact solution to the quantum Harmonic oscillator
problem, see Eq. (\ref{harmonic}), whereas, the Fourier method can only provide approximations to this 
solution. Therefore, rigorous numerical analysis should always be done case by case, with detailed 
consideration of functional class and various constraints.

\item {\sc [``Proof''- Section 2]}
``In practice, the requirement that $u(x)$ decays exponentially for $|x| \rightarrow \infty$ 
implies that the grid can be truncated to some large-but-finite span; the derivative sums 
are then truncated to summations over the entire truncated grid.'' 

{\sc [Analysis]}
The Fourier basis function $u(x)= e^{ikx}$ is widely  used in practice as well as in ``Proof'', 
and does not  decays exponentially for $|x| \rightarrow \infty$.

\item {\sc [``Proof''- Section 3]}
``However, the problem of summing slowly convergent series is an ancient one. A broad collection 
of schemes, known variously as ``summability'', ``sequence acceleration'' or ``sum-acceleration'' 
methods have been developed. Boyd [5] was the first to apply such ideas to pseudospectral 
series to invent the form of nonstandard differences called ``sum-accelerated pseudospectral''.''

{\sc [Analysis]}
The general idea of accelerating the Whittaker-Shannon-Kotel'nikov sampling, i.e., 
the sinc pseudospectral series, with a weight function was introduced as early as 1919 by M. 
Theis \cite{Theis}. Campbell introduced a weight of the form 
$w(x)=\frac{\int_{-1}^1\exp[(1-t^2)^{-1}-itx]dt}{\int_{-1}^1\exp[(1-t^2)^{-1}-itx]dt}$ to 
the sinc pseudospectral series in 1968 \cite{Campbell}.

\item {\sc [``Proof''- Section 3]}
``The sinc differentiation eigenvalue for the first derivative is the usual sine series $\cdots$.
The first option to obtain a sparse stencil is to simply {\it truncate} the infinite 
series at some upper limit $n$. This is a really bad idea because the error in the series, 
truncated after $n$ terms, is $O(1/n)$, and therefore unacceptably large.''

{\sc [Analysis]}
Note that ``the sinc differentiation eigenvalue for the first derivative'' refers to 
the result of differentiating $e^{(ikx)}$ with the Sinc method.  The truncation error of the 
sinc pseudospectral method depends not only on the decay rate of the differentiation kernel, 
but also on the decay rate of the function being differentiated. This partially 
explains why the performance in differentiating $e^{(ikx)}$ does not translate directly 
into the performance in differentiating other functions. In Section \ref{quantum} we 
show that for a quantum eigenvalue problem whose solution is of  Schwartz class, the sinc 
collocation outperforms the Boyd's Euler-accelerated sinc (Euler) and spectrally 
weighted least square differences (Sech) by many orders of magnitude, see 
Tables \ref{table.ho50} and \ref{table.ho200}.

\item {\sc [``Proof''- Section 5]}
``The shaded region shows that when $a < a_{\rm FD}$, the finite difference error is worse than that 
for the Gaussian weighting – but only for $K > \pi/2$.'' 

{\sc [Analysis]}
Section \ref{deriva} provides a counterexample. 
Ironically,  many more counterexamples appear in  Fig. 2 of ``Proof''.

\item {\sc [``Proof''- Section 5]}
``Given that the Fourier spectra of smooth functions fall off exponentially as $|K|$ increases, 
even small errors near $K = 0$ are intolerable.''

{\sc [Analysis]}
There is a simple counterexample: The Fourier spectra of smooth functions, such as 
the $C^{\infty}$ function of $e^{-x}\sin(x)$, do not fall off exponentially as $|K|$ 
increases.

\item {\sc [``Proof''- Section 5]}
``Why are finite differences so superior to DSC, given that both employ weights that are 
close-to-Gaussian? The answer is that the finite difference weights are tuned to give 
maximum accuracy in the limit $K \rightarrow 0$.''

{\sc [Analysis]}
It is well-known that high order 
finite difference weights are {\it exact} in  differentiating polynomials of the same order or
less, which have a wide Fourier spectral distribution beyond the 
component of $K \rightarrow 0$.

\item {\sc [``Proof''- Section 7]}
``How can the DSC/LDAF algorithms be salvaged? One could of course try a 
different sum-acceleration weighting from the Gaussian. However, this is merely to explore 
various instances of the sum-acceleration methods of [5]. ''

{\sc [Analysis]}

The many examples given here give the reader sufficient information 
to judge the validity of this claim.   

\end{itemize}

\section{Concluding remarks}

This paper examines the validity of claims made in Boyd's paper, ``A proof that the discrete 
singular convolution (DSC)/Lagrange-distributed approximation function (LDAF) method is inferior
to high order finite differences'' \cite{boyd06}, which is referred to as ``Proof''.  A wide variety 
of test problems are employed to compare six different numerical methods, including the 
discrete singular convolution with the regularized Shannon kernel (DSC-RSK), the  
standard finite difference (FD),
Boyd's  spectrally-weighted difference with the sech weight (Sech) \cite{Boyd94}, 
Boyd's finite difference (Boyd's FD) \cite{Boyd94},
Boyd's Euler-accelerated sinc algorithm (Euler) \cite{Boyd91,Boyd94}, and
Boyd's modified Euler-accelerated sinc algorithm (MEuler) \cite{Boyd91,Boyd94}.
These methods are employed in our comparisons because ``Proof''  has placed great emphasis 
of them, including many detailed expressions and/or numerical procedures, and because
they are related to some of ``Proof'''s claims. 
Some of the test problems, such as differentiating $e^{ikx}$ and a boundary value problem, 
are the designated tests used in either ``Proof''  or Boyd's earlier literature \cite{Boyd94}. 
In numerous examples, the DSC-RSK method outperforms all the other five methods by  factors of 
multiple orders of magnitude. In particular, for many problems with large wavenumbers or a 
wide range of wavenumbers (including large amplitudes in small wavenumbers), 
the DSC-RSK method is up to a million times more accurate than any of the other five methods mentioned. 

Although ``Proof'' presents error expressions for both the FD and DSC-RSK methods for 
differentiating the Fourier basis function, its claims are not directly supported by 
these error expressions, and the claim   about the superiority of the FD method over 
the DSC-RSK method is based on numerical experiments with  a limited set of parameters, 
and on some informal arguments. Moreover, ``Proof'''s claim about the superiority of 
spectrally-weighted differences over the DSC-RSK method is not supported by any analysis.

While we demonstrate that the DSC-RSK method outperforms the FD method for 
differentiating $e^{(ikx)}$  with fairly small wavenumbers, it may be the case 
that the FD method is superior with very small wavenumbers. 
In fact, the DSC-RSK method's less accurate performance in differentiating functions with 
{\it solely} very small wavenumbers is nothing new to us. It is for this reason that the DSC-RSK method was not 
proposed as another finite difference scheme.  Instead, it was proposed as 
a local spectral method, to be used for problems that are difficult for both 
low order methods and global spectral methods.  High frequency problems in flows, 
structural vibrations, and electromagnetic wave scattering and propagation are 
typical examples \cite{BaWeZh04,Lim,Weijsv,Zhaoijss,Zhou,Zhouijnmf}.
The DSC algorithm had hardly been used in small stencils in its applications, 
except for a couple of cases in complying with referee's requests. Nevertheless, 
for certain problems, it could outperform the FD method with small stencils as shown in Section \ref{deriva} 
and many other examples. The analysis of spectral convergence of the DSC-RSK method was given for 
certain class of functions in \cite{Qian}.

One may expect that for some physical problems where functions have  Fourier spectra 
with amplitudes decaying exponentially,
the ability to accurately approximate large wavenumbers is not important.  
In Section \ref{expon}, we present  an example in which the function has 
exponentially decaying amplitudes in wavenumbers $k$, i.e., amplitudes in large wavenumbers
are exponentially small comparing to those in small wavenumbers. Contrary to ``Proof'''s
claims, the DSC-RSK  method outperforms the FD method over all the stencils examined, and is up to 
$10^6$ times more accurate than the FD method and Boyd's spectrally weighted difference, the Sech method,
with some large stencils. This example indicates that for six numerical methods examined,
the dominant errors in derivatives originate from large wavenumbers. Therefore, a
method that is not only accurate for low wavenumbers, but also able to deliver high 
accuracy for large wavenumbers, will be more useful for many physical problems whose
Fourier spectra decay exponentially.

Indeed, many problems examined in this work involve considerable amplitudes in either  
large wavenumbers, or a wide range of wavenumbers which include both small and large wavenumbers.  
These problems are not suitable  for low order methods -- they require high order methods with relatively 
large stencils to achieve highly accurate results. In designing counterexamples, emphasis was not given to 
physical origins, partially because the claims of ``Proof''   have very little to do 
with applications, and partially because the DSC-RSK method has been applied to many 
practical problems in the past, in particular, to problems involving large wavenumbers
\cite{BaWeZh04,Lim,Weijsv,Zhaoijss,Zhou,Zhouijnmf}. The importance of large wavenumbers, 
or `short waves', cannot be overemphasized in scientific and engineering applications. 
At the beginning of the new millennium, Zienkiewicz \cite{Zienkiewicz} listed the problem 
of short waves in acoustics, electromagnetics or surface wave applications as one of two unsolved 
computational problems. Babu{\v s}ka and coworkers \cite{Babuska,Babuska2,Deraemaeker,Ihlenburg1,Ihlenburg2}
devoted much effort to constructing advanced numerical methods for high frequency waves in 
the Helmholtz equation.  Engquist and his co-workers \cite{Engquist} proposed  the segment projection
method for the propagation of high-frequency waves  in  waveguides.
Shu and Osher \cite{Shu} designed highly oscillatory problems in 
hyperbolic conservation laws to validate high-order shock-capturing methods. 
In space science and aerospace engineering, it is pertinent to quote Langley and Bardell \cite{Langley}
from their review  paper: ``...the prediction of medium to high frequency vibration levels 
is a particularly difficult task. ...there is no single technique 
which can be applied with confidence to all types of aerospace 
structures. Furthermore, there are certain problems of pressing 
practical concern for which it is not possible at present to make a 
reliable design prediction of high frequency vibration levels 
--- the prediction of on-orbit micro-vibration levels in satellite 
structures is arguably a problem of this type''.

Like Boyd's Sech method, the DSC-RSK method has a free parameter $r$. Since ``Proof'''s   
claim was very strong, specifically including all $a < a_{\rm  FD}$ (i.e., all $r > 
\frac{1} {\sqrt{2} a_{\rm  FD}}$, see ``Proof'''s   Abstract), we have chosen near optimal $r$ 
values in this work. In practice, $r$ values can be optimized according to specific 
applications. For a given stencil, there is a quite wide range of $r$ values that deliver very good 
results. Hence, the DSC-RSK method is relatively robust. Although it takes some experience and 
understanding to choose a near optimal $r$ value, there is no need to know the exact solution a prior. 
If one wishes to obtain a near optimal $r$ value, one can 
analyze the Fourier frequency response of the numerical solution. A near optimal $r$ value can 
then be obtained by  elevating (or  decreasing) the first choice of $r$ if the Fourier frequency 
response involves very high frequency components (or solely low frequency components).  

The sole purpose of this paper is to analyze the validity of ``Proof''. 
Although the reported counterexamples are a very small fraction of counterexamples we know,
we have no intention to claim the superiority of the DSC-RSK algorithm. Given the great 
diversity of problems with different physical origins, it is improper to claim
that one method is superior to others {\it in general} without detailed analysis and
comparison. However, for a given problem, one could show that some methods are more suitable 
than others. Case-by-case study is very important for validating new numerical methods.
In view of  the fact that detailed comparisons between the DSC-RSK method  and  many other 
numerical methods, including many other spectrally-weighted  differences,
have not been made anywhere, it is entirely 
possible to find another method that outperforms the DSC-RSK method for some examples studied in this work. 
It is also possible to come up with other examples for which the DSC-RSK method does not perform as well as the 
other five finite difference-type  methods examined in this work. However, these possibilities 
do not affect the conclusion of the present paper, which demonstrates that the general statements 
in ``Proof'' are unfounded.

\vspace*{1.5cm}

\centerline{\bf Acknowledgments}

{This work was supported in part by NSF Grant IIS-0430987. 
}

%%%%%%%%%%%%%%%%%%%%%%%%%%%%% 
%%%%%%%%%%%%%%%%%%%%%%%%%%%%% 
%%% Start of Bibliography %%% 
%%%%%%%%%%%%%%%%%%%%%%%%%%%%% 
 
%%%%%%%%%%%%%%%%%%%%%%%%%%% 
%%% End of Bibliography %%% 
%%%%%%%%%%%%%%%%%%%%%%%%%%% 

%\end{document}

\end{document}